\documentclass[
submission
]{dmtcs-episciences}

\usepackage[utf8]{inputenc}
\usepackage{subfigure}
\usepackage[round]{natbib}
\newtheorem{thm}{Theorem}[section]
\newtheorem{cor}[thm]{Corollary}
\newtheorem{lem}[thm]{Lemma}

\newtheorem{obs}[thm]{Observation}
\newtheorem{pro}[thm]{Proposition}

\usepackage[round]{natbib}

\author{Yaping Mao\affiliationmark{1,4}\thanks{Supported by the National Science Foundation of China
(Nos. 11601254, 11551001, 11161037, and 11461054) and the Science
Found of Qinghai Province (Nos. 2016-ZJ-948Q, and 2014-ZJ-907).}
  \and Eddie Cheng\affiliationmark{2,4}
  \and Zhao Wang\affiliationmark{3}\thanks{Corresponding author}}
\title[Steiner distance in product networks]{Steiner distance in product networks}
\affiliation{
  School of Mathematics and Statistics, Qinghai Normal University, Xining, Qinghai 810008, China\\
  Department of Mathematics and Statistics, Oakland University, Rochester, MI USA 48309\\
  College of Science, China Jiliang University, Hangzhou 310018, China\\
  Center for Mathematics and Interdisciplinary Sciences of Qinghai Province, Xining, Qinghai 810008, China}
\received{2017-3-8} \revised{2018-1-12, 2018-7-10, 2018-9-19}
\accepted{2018-9-19}
\begin{document}
\publicationdetails{20}{2018}{2}{8}{3175}
\maketitle
\begin{abstract}
For a connected graph $G$ of order at least $2$ and $S\subseteq
V(G)$, the \emph{Steiner distance} $d_G(S)$ among the vertices of
$S$ is the minimum size among all connected subgraphs whose vertex
sets contain $S$. Let $n$ and $k$ be two integers with $2\leq k\leq
n$. Then the \emph{Steiner $k$-eccentricity $e_k(v)$} of a vertex
$v$ of $G$ is defined by $e_k(v)=\max \{d_G(S)\,|\,S\subseteq V(G),
\ |S|=k, \ and \ v\in S \}$. Furthermore, the \emph{Steiner
$k$-diameter} of $G$ is $sdiam_k(G)=\max \{e_k(v)\,|\, v\in V(G)\}$.
In this paper, we investigate the Steiner distance and Steiner
$k$-diameter of Cartesian and lexicographical product graphs. Also,
we study the Steiner $k$-diameter of some networks.
\end{abstract}

\section{Introduction}
\label{sec:in} In this paper, we consider graphs that are
undirected, finite and simple. We refer the readers to \cite{Bondy}
for graph theoretical notations and terminology that are not defined
here. For a graph $G$, let $V(G)$, $E(G)$, and $\delta(G)$ denote
the set of vertices, the set of edges and the minimum degree of $G$,
respectively. We refer to $|V(G)|$ the order of the graph and
$|E(G)|$ the size of the graph. The degree of a vertex $v$ in $G$ is
denoted by $deg_G(v)$. In this paper, $K_{n}$, $P_n$, $K_{1,n-1}$
and $C_n$ correspond to the complete graph of order $n$, the path of
order $n$, the star of order $n$, and the cycle of order $n$,
respectively. If $X\subseteq V(G)$, we use $G[X]$ to denote the
subgraph induced by $X$. Similarly, if $F\subseteq E(G)$, let $G[F]$
denote the subgraph induced by $F$. If $X\subseteq V(G)\cup E(G)$,
we use $G-X$ to denote the subgraph of $G$ obtained from $G$ by
removing all the elements of $X$ and the edges incident to vertices
that are in $X$. If $X=\{x\}$, we write $G-x$ for notational
simplicity. For $X,Y\subseteq V(G)$, we use $E_G[X,Y]$ to denote the
set of edges of $G$ with one end in $X$ and the other end in $Y$. If
$X=\{x\}$, we simply write $E_G[x,Y]$ for $E_G[\{x\},Y]$. We divide
our introduction into subsections to state the motivations of this
paper.

\subsection{Distance and its generalizations}

Distance is a fundamental concept in graph theory. Let $G$ be a
connected graph. The \emph{distance} between two vertices $u$ and
$v$ in $G$ is the length of a shortest path between them, and it is
denoted by $d_G(u,v)$. The \emph{eccentricity} of $v$ in $G$,
denoted by $e_G(v)$ (or simply $e(v)$ if it is clear from the
context), is $\max\{d_G(u,v)\,|\,u\in V(G)\}$. In addition, we
define the \emph{radius} $rad(G)$ and the \emph{diameter} $diam(G)$
of $G$ to be $rad(G)=\min\{e(v)\,|\,v\in V(G)\}$ and $diam(G)=\max
\{e(v)\,|\,v\in V(G)\}$. It is a standard exercise to check that
$rad(G)\leq diam(G) \leq 2 rad(G)$. The \emph{center} $C(G)$ of $G$
is the subgraph induced by the vertices with eccentricity equal to
the radius. For more details on distance, we refer to \cite{Harary,
Goddard}.

We observe that the distance between two vertices $u$ and $v$ in $G$
is equal to the minimum size of a connected subgraph of $G$
containing both $u$ and $v$. This suggests a generalization of the
concept of distance. The Steiner distance of a graph, introduced by
\cite{Chartrand} in 1989, is
such a natural and nice generalization. Let $S$ be a set of vertices
in a graph $G(V,E)$ where $|S|\geq 2$. We define \emph{an
$S$-Steiner tree} or \emph{a Steiner tree connecting $S$} (or
simply, \emph{an $S$-tree}) to be a subgraph $T(V',E')$ of $G$ that
is a tree with $S\subseteq V'$. Moreover, the \emph{Steiner
distance} $d_G(S)$ of $S$ in $G$ (or simply the distance of $S$) is
the minimum size among all connected subgraphs whose vertex sets
contain $S$. (Set $d_G(S)=\infty$ when there is no $S$-Steiner tree
in $G$.) We remark that if $H$ is a connected subgraph of $G$ such
that $S\subseteq V(H)$ and $|E(H)|=d_G(S)$, then $H$ is a tree. We
further remark that $d_G(S)=\min\{e(T)\,|\,S\subseteq V(T)\}$, where
$T$ is subtree of $G$. Finally, if $S=\{u,v\}$, then $d_G(S)=d(u,v)$
is the classical distance between $u$ and $v$. The following
observation is obvious.
\begin{obs}\label{obs1-1}
Let $G$ be a graph of order $n$ and $k$ be an integer with $2\leq
k\leq n$. If $S\subseteq V(G)$ and $|S|=k$, then $d_G(S)\geq k-1$.
\end{obs}

Let $n$ and $k$ be two integers with $2\leq k\leq n$. We define the
\emph{Steiner $k$-eccentricity $e_k(v)$} of a vertex $v$ of $G$ to
be $e_k(v)=\max \{d(S)\,|\,S\subseteq V(G), |S|=k,~and~v\in S \}$,
the \emph{Steiner $k$-radius} of $G$ to be $srad_k(G)=\min \{
e_k(v)\,|\,v\in V(G)\}$, and the \emph{Steiner $k$-diameter} of $G$
is $sdiam_k(G)=\max \{e_k(v)\,|\,v\in V(G)\}$. We remark that for
every connected graph $G$ that $e_2(v)=e(v)$ for all vertices $v$ of
$G$ and that $srad_2(G)=rad(G)$ and $sdiam_2(G)=diam(G)$. It is not
difficult to see the following observation.

\begin{obs}\label{obs1-2}
Let $k,n$ be two integers with $2\leq k\leq n$.

$(1)$ If $H$ is a spanning subgraph of $G$, then $sdiam_k(G)\leq
sdiam_k(H)$.

$(2)$ For a connected graph $G$, $sdiam_k(G)\leq sdiam_{k+1}(G)$.
\end{obs}

\cite{ChartrandOZ} obtained the following
upper and lower bounds of $sdiam_k(G)$.
\begin{thm}{\upshape\cite{ChartrandOZ}}\label{th1-3}
Let $k,n$ be two integers with $2\leq k\leq n$, and let $G$ be a
connected graph of order $n$. Then $k-1\leq sdiam_k(G)\leq n-1$.
Moreover, the upper and lower bounds are sharp.
\end{thm}

\cite{DankelmannSO2} showed that
$sdiam_k(G)\leq \frac{3|V(G)|}{\delta(G)+1}+3k$. \cite{AliDM} improved the bound and showed that
$sdiam_k(G)\leq \frac{3|V(G)|}{\delta(G)+1}+2k-5$ where $G$ is
connected. Moreover, they showed that these bounds are
asymptotically best possible via a construction.

\subsection{Related concepts}

Although we will not consider these related concepts in this paper,
they provide a context of problems related to Steiner distance. As a
generalization of the center of a graph, one defines the
\emph{Steiner $k$-center} $C_k(G)\ (k\geq 2)$ of a connected graph
$G$ to be the subgraph induced by the vertices $v$ of $G$ where
$e_k(v)=srad_k(G)$. \cite{OellermannT} showed
that every graph is the $k$-center of some graph. Moreover, they
showed that the $k$-center of a tree is a tree and they
characterized those trees that are $k$-centers of trees. The
\emph{Steiner $k$-median} of $G$ is the subgraph of $G$ induced by
the vertices of $G$ of minimum Steiner $k$-distance. The papers
\cite{Oellermann, Oellermann2, OellermannT} contain important
results for Steiner centers and Steiner medians. For more details on
the Steiner distance parameters, we refer to the survey paper
\cite{MaoSurvey} and papers \cite{Ali, Caceresa, Moscarini,
DankelmannE, DankelmannSO2, DayOS, GoddardOS, MaoMC}.

Let $G$ be a $k$-connected graph and $u$, $v$ be a pair of vertices
of $G$. Let $P_k(u,v)=\{P_1,P_2,\cdots,P_k\}$ be a family of $k$
internally vertex-disjoint paths between $u$ and $v$ and
$l(P_k(u,v))$ be the length of the longest path in $P_k(u,v)$. Then
the \emph{$k$-distance} $d_k(u,v)$ between vertices $u$ and $v$ is
the smallest $l(P_k(u,v))$ among all $P_k(u,v)$'s and the
\emph{$k$-diameter} $d_k(G)$ of $G$ is the maximum $k$-distance
$d_k(u,v)$ over all pairs $u,v$ of vertices of $G$. The concept of
$k$-diameter has its origin in the analysis of routings in networks
as described by \cite{Chung,Du,Hsu,Hsu2,Meyer}.

Perhaps the most famous Steiner type problem is the Steiner tree
problem. The original Steiner tree problem was stated for the
Euclidean plane: Given a set of points on the plane, the goal is to
connect these points, and possibly additional points, by line
segments between some pairs of these points such that the total
length of these line segments is minimized. The graph theoretical
version~\cite{Hakimi,Levi} is as follows: Given a graph and a set of
vertices $S$, find a connected subgraph with minimum number of edges
that contains $S$. This is, in general, an NP-hard problem
\cite{HwangRW}. There is also a corresponding weighted version.
Obviously, this has applications in computer science and  electrical
engineering. For example, a graph can be a computer network with
vertices being computers and edges being links between them. Here
the Steiner tree problem is to find a subnetwork containing these
computers with the least number of links. We can replace processors
by electrical stations for applications in electrical networks.

\cite{LMG} gave such a concept. They defined the {\it
$k$-center Steiner Wiener index\/} $SW_k(G)$ of the graph $G$ to be
$$
SW_k(G)=\sum_{S\subseteq V(G), \ |S|=k} d(S)\,.
$$
For $k=2$, it coincides with the ordinary Wiener index. One usually
considers $SW_k$ for $2 \leq k \leq n-1$. However, the above
definition can be extended to $k=1$ and $k=n$ as well where
$SW_1(G)=0$ and $SW_n(G)=n-1$. There are other related concepts such
as the Steiner Harary index.  Both indices have chemical
applications~\cite{FurtulaGK, GFL}. In addition,
\cite{GutmanSDD} gave a generalization of the concept of degree
distance, and then \cite{MaoDas} gave a generalization
of the concept of Gutman index. We refer the readers to
\cite{FurtulaGK, GFL, GutmanSDD, LMG, LMG2, MaoDas, MWG, MWGK, MWGL}
for details.

\subsection{Products of graphs}
The main focus of this paper is Steiner $k$-diameter of two products
of graphs, namely, the Cartesian product and the lexicographic
product. These are well-known products. See \cite{Hammack}.

$\bullet$ The {\it Cartesian product} of two graphs $G$ and $H$,
written as $G\Box H$, is the graph with vertex set $V(G)\times
V(H)$, in which two vertices $(g,h)$ and $(g',h')$ are adjacent if
and only if $g=g'$ and $(h,h')\in E(H)$, or $h=h'$ and $(g,g')\in
E(G)$.

$\bullet$ The {\it lexicographic product} of two graphs $G$ and $H$,
written as $G\circ H$, is defined as follows: $V(G\circ
H)=V(G)\times V(H)$, and two distinct vertices $(g,h)$ and $(g',h')$
of $G\circ H$ are adjacent if and only if either $(g,g')\in E(G)$ or
$g=g'$ and $(h,h')\in E(H)$.

It is easy to see that the Cartesian product is commutative, that
is, $G\Box H$ is isomorphic to $H\Box G$. However, the lexicographic
product is non-commutative.

Product networks are important as often the resulting graph inherits
properties from its factors. Both the lexicographical product and
the Cartesian product are important concepts. See \cite{Bao, DayA,
Hammack, Ku}.

\cite{Gologranc} obtained a sharp lower bound for
Steiner distance of Cartesian product graphs. We continue this study
in Section $2$ by obtaining a sharp upper bound for Steiner
distance. In addition, we will also present sharp upper and lower
bounds for Steiner $k$-diameter of Cartesian product graphs. In
Section $3$, we derive the results for Steiner distance and Steiner
$k$-diameter of lexicographic product graphs, which strengthen a
result given by \cite{AnandCKP}. In Section $4$, we
give some applications of our main results, and study the Steiner
diameter of some important networks.

\section{Results for Cartesian product}

In this paper, let $G$ and $H$ be two graphs with
$V(G)=\{g_1,g_2,\ldots,g_{n}\}$ and $V(H)=\{h_1,h_2,\ldots,h_{m}\}$,
respectively. Then $V(G\ast H)=\{(g_i,h_j)\,|\,1\leq i\leq n, \
1\leq j\leq m\}$, where $\ast$ denotes the Cartesian product
operation or lexicographical product operation. For $h\in V(H)$, we
use $G(h)$ to denote the subgraph of $G\ast H$ induced by the vertex
set $\{(g_i,h)\,|\,1\leq i\leq n\}$. Similarly, for $g\in V(G)$, we
use $H(g)$ to denote the subgraph of $G\circ H$ induced by the
vertex set $\{(g,h_j)\,|\,1\leq j\leq m\}$.

The following observation can be easily seen.
\begin{obs}\label{obs2-1}
Let $G$ be a connected graph, and let $S\subseteq V(G)$ and $|S|=3$.
Let $T$ be a minimal $S$-Steiner tree in $G$. Then the tree $T$
satisfies one of the following conditions.

$\bullet$ $T$ is a path;

$\bullet$ $T$ is a subdivision of $K_{1,3}$.
\end{obs}

We start with the following basic result.
\begin{lem}{\upshape \cite{Hammack}}\label{lem2-1}
Let $G$ and $H$ be two graphs, and let $(g,h)$ and $(g',h')$ be two
vertices of $G\Box H$. Then
$$
d_{G\Box H}((g,h),(g',h'))= d_{G}(g,g')+ d_{H}(h,h').
$$
\end{lem}

\subsection{Steiner distance of Cartesian product graphs}

\cite{Gologranc} obtained the following lower bound
for Steiner distance.
\begin{lem}{\upshape \cite{Gologranc}}\label{lem2-2}
Let $k\geq 2$ be an integer, and let $G,H$ be two connected graphs.
Let $S=\{(g_{i_1},h_{j_1}),$
$(g_{i_2},h_{j_2}),\ldots,(g_{i_k},h_{j_k})\}$ be a set of distinct
vertices of $G\Box H$. Let $S_G=\{g_{i_1},g_{i_2},\ldots,g_{i_k}\}$
and $S_H=\{h_{j_1},h_{j_2},\ldots,h_{j_k}\}$. Then
$$
d_{G\Box H}(S)\geq d_G(S_G)+d_H(S_H).
$$
\end{lem}

We will show that the inequality in Lemma \ref{lem2-2} can be
equality if $k=3$; shown in following Corollary \ref{cor2-2}. But,
for general $k \ (k\geq 4)$, from Lemma \ref{lem2-2} and Corollary
\ref{cor2-2}, one may conjecture that for two connected graphs
$G,H$, $d_{G\Box H}(S)=d_G(S_G)+d_H(S_H)$, where
$S=\{(g_{i_1},h_{j_1}),(g_{i_2},h_{j_2}),\ldots,(g_{i_k},h_{j_k})\}\subseteq
V(G\Box H)$, $S_G=\{g_{i_1},g_{i_2},\ldots,g_{i_k}\}\subseteq V(G)$
and $S_H=\{h_{j_1},h_{j_2},\ldots,h_{j_k}\}\subseteq V(H)$.\\

\noindent{\bf Remark 1:} Actually, the equality $d_{G\Box
H}(S)=d_G(S_G)+d_H(S_H)$ is not true for $|S|\geq 4$. For example,
let $G$ be a tree with degree sequence $(3,2,1,1,1)$ and $H$ be a
path of order $5$. Let
$S=\{(g_1,h_1),(g_2,h_2),(g_3,h_3),(g_4,h_4)\}$ be a vertex set of
$G\Box H$ shown in Fig.1. Then $d_G(S_G)=4$ for
$S_G=\{g_1,g_2,g_3,g_4\}$, and $d_H(S_H)=4$ for
$S_H=\{h_1,h_2,h_3,h_4\}$. One can check that there is no
$S$-Steiner tree of size $8$ in $G\Box H$, which implies $d_{G\Box
H}(S)\geq 9$.
\begin{figure}[htbp]
  \begin{center}
    \includegraphics[scale=0.7]{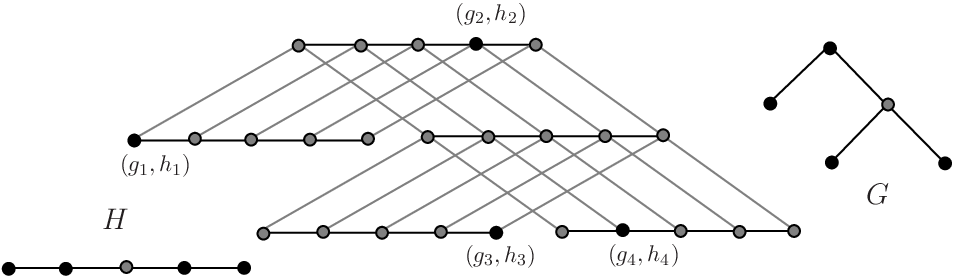}
    \caption{Graphs for Remark 1.}
    \label{fig:logo}
  \end{center}
\end{figure}

Although the conjecture of such an ideal formula is not correct, it
is possible to give a strong upper bound for general $k \ (k\geq
3)$. Remark 1 also indicates that obtaining a nice formula for the
general case may be difficult. We now give such an upper bound of
$d_{G\Box H}(S)$ for $S\subseteq V(G\Box H)$ and $|S|=k$.
\begin{thm}\label{th2-1}
Let $k,m,n$ be three integers with $3\leq k\leq mn$, and let $G,H$
be two connected graphs with $V(G)=\{g_1,g_2,\ldots,g_{n}\}$ and
$V(H)=\{h_1,h_2,\ldots,h_{m}\}$. Let
$S=\{(g_{i_1},h_{j_1}),(g_{i_2},h_{j_2}),\ldots,(g_{i_k},h_{j_k})\}$
be a set of distinct vertices of $G\Box H$, $S_G=\{g_{i_1},g_{i_2},$
$\ldots,g_{i_k}\}$, and $S_H=\{h_{j_1},h_{j_2},\ldots,h_{j_k}\}$,
where $S_G\subseteq V(G)$, $S_H\subseteq V(H)$ ($S_G,S_H$ are both
multi-sets). Then
\begin{eqnarray*}
d_G(S_G)+d_H(S_H)&\leq &d_{G\Box H}(S)\\[0.1cm]
&\leq&\min\{d_G(S_G)+(r+1)d_H(S_H),d_H(S_H)+(t+1)d_G(S_G)\},
\end{eqnarray*}
where $r,t \ (0\leq r,t\leq k-3)$ are defined as follows.

$\bullet$ Let $X_G^{i} \ (1\leq i\leq {k\choose 3})$ be all the
$(k-3)$-multi-subsets of $\{g_{i_1},g_{i_2},,\ldots,g_{i_k}\}$ in
$G$, and let $r_i$ be the numbers of distinct vertices in $X_G^{i} \
(1\leq i\leq {k\choose 3})$, and let $r=\min\{r_i\,|\,1\leq i\leq
{k\choose 3}\}$.

$\bullet$ Let $Y_H^{j} \ (1\leq j\leq {k\choose 3})$ be all the
$(k-3)$-multi-subsets of $\{h_{j_1},h_{j_2},\ldots,h_{j_k}\}$ in
$H$, and let $t_j$ be the numbers of distinct vertices in $Y_H^{j} \
(1\leq j\leq {k\choose 3})$, and let $t=\min\{t_j\,|\,1\leq j\leq
{k\choose 3}\}$.
\end{thm}
\begin{pf}
From Lemma \ref{lem2-2}, we have $d_{G\Box H}(S)\geq
d_G(S_G)+d_H(S_H)$. By symmetry, we only need to show $d_{G\Box
H}(S)\leq d_G(S_G)+(r+1)d_H(S_H)$. Recall that
$V(G)=\{g_1,g_2,\ldots,g_{n}\}$ and $V(H)=\{h_1,h_2,\ldots,h_{m}\}$.
Without loss of generality, we assume that
$H(g_{1}),H(g_{2}),\ldots,$ $H(g_{a})$ be the $H$ copies such that
$|V(H(g_{i}))\cap S|\neq 0$, $1\leq i\leq a$. Then
$(g_{i_{1}},h_{j_{1}}),(g_{i_{2}},h_{j_{2}}),$
$\ldots,(g_{i_{k}},h_{j_{k}})\in \bigcup_{i=1}^aV(H(g_i))$, and
hence we have the following cases to consider.

{\bf Case 1.} For each $H(g_{i}) \ (1\leq i\leq a)$,
$|V(H(g_{i}))\cap S|\geq 2$.

Without loss of generality, let $V(H(g_{1}))\cap
S=\{(g_{i_1},h_{j_1}),(g_{i_2},h_{j_2}),\ldots,(g_{i_s},h_{j_s})\}$,
where $s\geq 2$. Thus, we have $(g_{i_p},h_{j_p})=(g_1,h_{j_p})$ for
each $p \ (1\leq p\leq s)$, and
$(g_{i_{s+1}},h_{j_{s+1}}),(g_{i_{s+2}},h_{j_{s+2}}),\ldots,$
$(g_{i_{k}},h_{j_{k}})\in \bigcup_{i=2}^aV(H(g_i))$. Note that
$(g_{1},h_{j_1}),(g_1,h_{j_2}),\ldots,(g_1,h_{j_s})\in V(H(g_{1}))$.
On one hand, since there is an $S_H$-Steiner tree of size $d_H(S_H)$
in $H$, it follows that there exists an Steiner tree of size
$d_H(S_H)$ connecting
\begin{eqnarray*}
&&\{(g_{1},h_{j_1}),(g_1,h_{j_2}),\ldots,(g_1,h_{j_s})\}\cup \{(g_1,h_{j_{s+1}}),(g_1,h_{j_{s+2}}),\ldots,(g_1,h_{j_k})\}\\[0.1cm]
&=&\{(g_{i_1},h_{j_1}),(g_{i_2},h_{j_2}),\ldots,(g_{i_s},h_{j_s})\}\cup
\{(g_1,h_{j_{s+1}}),(g_1,h_{j_{s+2}}),\ldots,(g_1,h_{j_k})\}
\end{eqnarray*}
in $H(g_1)$, say $T(g_1)$. For each $i \ (2\leq i\leq k)$, let
$T(g_i)$ be the Steiner tree in $H(g_i)$ corresponding to $T(g_1)$
in $H(g_1)$. Note that $T(g_i) \ (1\leq i\leq k)$ is the Steiner
tree of size $d_H(S_H)$ connecting
$\{(g_{i},h_{j_1}),(g_i,h_{j_2}),\ldots,(g_i,h_{j_s}),(g_i,h_{j_{s+1}}),(g_i,h_{j_{s+2}}),\ldots,(g_i,h_{j_k})\}$
in $H(g_i)$. One can see that
$(g_{i_{s+1}},h_{j_{s+1}}),\ldots,(g_{i_{k}},h_{j_{k}})\in
\bigcup_{i=2}^aV(T(g_i))$. On the other hand, since there is an
$S_G$-Steiner tree of size $d_G(S_G)$ in $G$, it follows that there
exists an Steiner tree of size $d_G(S_G)$ connecting
$\{(g_1,h_{j_1}),$ $(g_2,h_{j_1}),\ldots,(g_a,h_{j_1})\}$ in
$G(h_{j_1})$, say $T(h_{j_1})$. Furthermore, the subgraph induced by
the edges in $\left(\bigcup_{i=1}^{a}E(T(g_i))\right)\cup
E(T(h_{j_1}))$ is an $S$-Steiner tree in $G\Box H$ (see Fig.$2$
$(a)$), and hence $d_{G\Box H}(S)\leq d_G(S_G)+ad_H(S_H)$.

From the definition of $r$, if $|V(H(g_{i}))\cap S|\geq 4$ for each
$H(g_{i}) \ (1\leq i\leq a)$, then $r=a$ and $d_{G\Box H}(S)\leq
d_G(S_G)+rd_H(S_H)$. If there exists some $H(g_{i}) \ (1\leq i\leq
a)$ such that $2\leq |V(H(g_{i}))\cap S|\leq 3$ for $H(g_{i}) \
(1\leq i\leq a)$, then $r=a-1$ and $d_{G\Box H}(S)\leq
d_G(S_G)+(r+1)d_H(S_H)$.
\begin{figure}[htbp]
  \begin{center}
    \includegraphics[scale=0.7]{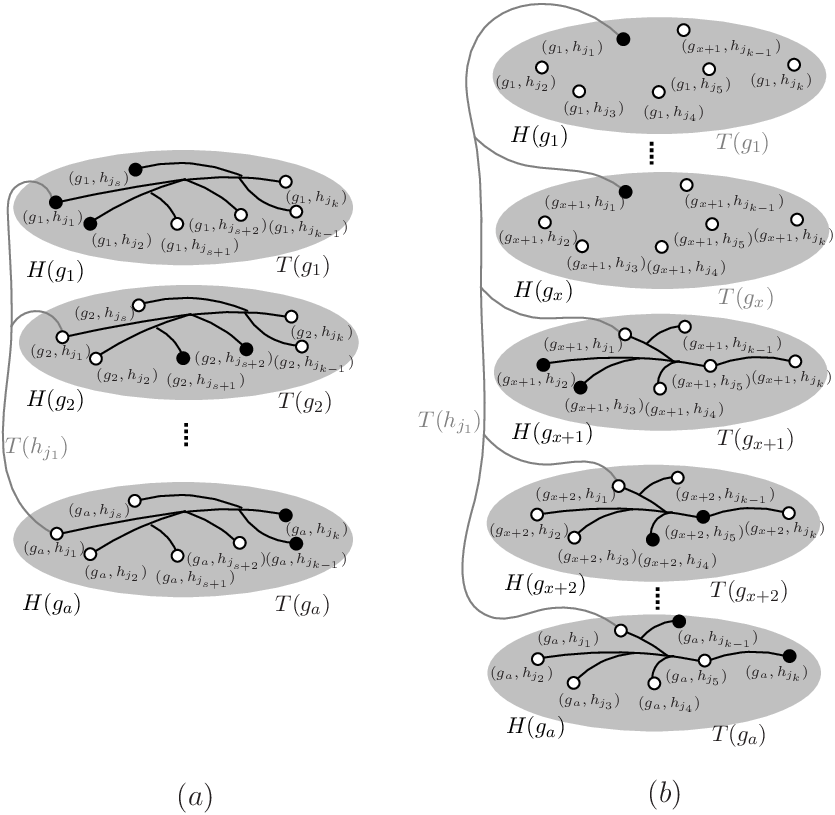}
    \caption{Graphs for Cases 1 and 2 in the proof of Theorem
\ref{th2-1}.}
    \label{fig:logo}
  \end{center}
\end{figure}

{\bf Case 2.} There exists some $H(g_{i})$ such that
$|V(H(g_{i}))\cap S|=1$, where $1\leq i\leq a$.

Without loss of generality, we assume that $|V(H(g_{i}))\cap S|=1$
for each $i \ (1\leq i\leq x)$, where $1\leq x\leq a$. For $x\neq
a$, we have $|V(H(g_{i}))\cap S|\geq 2$ for each $i \ (x+1\leq i\leq
a)$. One can see that
$$
x=|\{H(g_{i})\,|\,|V(H(g_{i}))\cap S|=1, \ 1\leq i\leq a\}|.
$$

{\bf Subcase 2.1.} $x\geq 3$.

If $|\{h_{j_1},h_{j_2},\cdots,h_{j_x}\}|=1$, then
$h_{j_1}=h_{j_2}=\cdots=h_{j_x}$. Since there is an $S_G$-Steiner
tree of size $d_G(S_G)$ in $G$, it follows that there exists an
Steiner tree of size $d_G(S_G)$ connecting
$\{(g_1,h_{j_1}),(g_2,h_{j_1}),$ $\ldots,(g_a,h_{j_1})\}$ in
$G(h_{j_1})$, say $T(h_{j_1})$. Since there is an $S_H$-Steiner tree
of size $d_H(S_H)$ in $H$, it follows that there exists an Steiner
tree of size $d_H(S_H)$ connecting $\{ (g_{x+1}, h_{j_{1}})\} \cup
\{(g_{x+1}, h_{j_{x+1}}),(g_{x+1}, h_{j_{x+2}}),$ $  \ldots,
(g_{x+1}, h_{j_{k}}) \}$ in $H(g_{x+1})$, say $T(g_{x+1})$. For each
$i \ (x+2\leq i\leq a)$, let $T(g_i)$ be the Steiner tree in
$H(g_i)$ corresponding to $T(g_{x+1})$ in $H(g_{x+1})$. Note that
$T(g_i) \ (x+1\leq i\leq a)$ is the Steiner tree of size $d_H(S_H)$
connecting
$\{(g_{i},h_{j_{x+1}}),(g_i,h_{j_{x+2}}),\ldots,(g_i,h_{j_k})\}$ in
$H(g_i)$. Furthermore, the subgraph induced by the edges in
$\left(\bigcup_{i=x+1}^{a}E(T(g_i))\right)\cup E(T(h_{j_1}))$ is an
$S$-Steiner tree (see Fig.2 $(b)$), and hence $d_{G\Box H}(S)\leq
d_G(S_G)+(a-x)d_H(S_H)\leq d_G(S_G)+(a-3)d_H(S_H)$. From the
definition of $r$, we have $r=a-3$, and hence $d_{G\Box H}(S)\leq
d_G(S_G)+rd_H(S_H)\leq d_G(S_G)+(r+1)d_H(S_H)$, as desired.

If $|\{h_{j_1},h_{j_2},\cdots,h_{j_x}\}|=2$, then we can assume that
$h_{j_1}=h_{j_2}=\ldots=h_{j_s}$,
$h_{j_{s+1}}=h_{j_{s+2}}=\ldots=h_{j_x}$, and $h_{j_1}\neq h_{j_x}$.
Furthermore, we can assume that $s\geq 2$. Since there is an
$S_H$-Steiner tree of size $d_H(S_H)$ in $H$, it follows that there
is a Steiner tree of size $d_H(S_H)$ connecting
$\{(g_{s+1},h_{j_{s+1}}),$
$(g_{s+1},h_{j_{x+2}}),\ldots,(g_{s+1},h_{j_k})\}$ in $H(g_{s+1})$,
say $T(g_{s+1})$. For each $i \ (s+2\leq i\leq a)$, let $T(g_i)$ be
the Steiner tree in $H(g_i)$ corresponding to $T(g_{s+1})$ in
$H(g_{s+1})$. Since there is an $S_G$-Steiner tree of size
$d_G(S_G)$ in $G$, it follows that there exists an Steiner tree of
size $d_G(S_G)$ connecting $\{(g_1,h_{j_1}),(g_2,h_{j_1}),$
$\ldots,(g_a,h_{j_1})\}$ in $G(h_{j_1})$, say $T(h_{j_1})$. Then the
subgraph induced by the edges in
$$
\left(\bigcup_{i=s+1}^{a}E(H(g_i))\right)\cup E(T(h_{j_1}))
$$
is an $S$-Steiner tree in $G\Box H$, and hence $d_{G\Box H}(S)\leq
d_G(S_G)+(a-s)d_H(S_H)\leq d_G(S_G)+(a-2)d_H(S_H)$. Since $r=a-3$,
it follows that $d_{G\Box H}(S)\leq
d_G(S_G)+(a-2)d_H(S_H)=d_G(S_G)+(r+1)d_H(S_H)$.

From now on, we assume $|\{h_{j_1},h_{j_2},\cdots,h_{j_x}\}|\geq 3$.
Note that there is an $S_H$-Steiner tree of size $d_H(S_H)$ in $H$,
say $T$. Without loss of generality, let $h_{j_1}\neq h_{j_2}\neq
h_{j_3}$. Since $h_{j_1},h_{j_2},h_{j_3}\in V(T)$, it follows that
there is a minimal subtree $T'$ connecting
$\{h_{j_1},h_{j_2},h_{j_3}\}$ in $T$. From Observation \ref{obs2-1},
$T'$ is a path or $T'$ is a subdivision of $K_{1,3}$. If $T'$ is a
path, then without loss of generality, we can assume $h_{j_2}$ is
the interval vertex of $T'$. Therefore, there are a unique
$(h_{j_1},h_{j_2})$-path, say $P^1$, and a unique
$(h_{j_2},h_{j_3})$-path, say $P^2$, in $T'$. If $T'$ is a
subdivision of $K_{1,3}$, then there exists a vertex in $T'$, say
$h^*\in V(H)\setminus \{h_{j_1},h_{j_2},h_{j_3}\}$, such that there
are three paths $Q^1,Q^2,Q^3$ connecting $h^*$ and
$h_{j_1},h_{j_2},h_{j_3}$, respectively, in $T'$.

We first consider the case that $T'$ is a path. On one hand, for
each $i \ (1\leq i\leq k)$, let $T(g_i)$ be the Steiner tree in
$H(g_i)$ corresponding to $T$ in $H$. Note that $T(g_i)$ is the
Steiner tree of size $d_H(S_H)$ connecting
$\{(g_{i},h_{j_1}),(g_i,h_{j_2}),\ldots,(g_i,h_{j_k})\}$ in
$H(g_i)$. For each $i \ (1\leq i\leq 3)$, let $P^1(g_i)$ be the path
in $H(g_i)$ corresponding to $P^1$ in $H$, and let $P^2(g_i)$ be the
path in $H(g_i)$ corresponding to $P^2$ in $H$. On the other hand,
since there is an $S_G$-Steiner tree of size $d_G(S_G)$ in $G$, it
follows that there exists an Steiner tree of size $d_G(S_G)$
connecting $\{(g_1,h_{j_2}),(g_2,h_{j_2}),$ $\ldots,(g_k,h_{j_2})\}$
in $G(h_{j_2})$, say $T(h_{j_2})$. Furthermore, the subgraph induced
by the edges in
$$
\left(\bigcup_{i=4}^{a}E(T(g_i))\right)\cup E(P^1(g_1))\cup
E(P^2(g_3))\cup E(T(h_{j_2}))
$$
is an $S$-Steiner tree in $G\Box H$ (see Fig.3 $(a)$), and hence
$d_{G\Box H}(S)\leq d_G(S_G)+(a-2)d_H(S_H)$. Since $r=a-3$, it
follows that $d_{G\Box H}(S)\leq d_G(S_G)+(r+1)d_H(S_H)$.

Next, we consider the case that $T'$ is a subdivision of $K_{1,3}$.
On one hand, for each $i \ (1\leq i\leq k)$, let $T(g_i)$ be the
tree in $H(g_i)$ corresponding to $T$ in $H$. Note that $T(g_i)$ is
the Steiner tree of size $d_H(S_H)$ connecting
$\{(g_{i},h_{j_1}),(g_i,h_{j_2}),\ldots,(g_i,h_{j_k})\}$ in
$H(g_i)$. For each $i \ (1\leq i\leq 3)$, let $Q^1(g_i)$ be the path
in $H(g_i)$ corresponding to $Q^1$ in $H$, and let $Q^2(g_i)$ be the
path in $H(g_i)$ corresponding to $Q^2$ in $H$, and let $Q^3(g_i)$
be the path in $H(g_i)$ corresponding to $Q^3$ in $H$. For each $i \
(1\leq i\leq k)$, let $(g_i,h^*)$ be the path in $H(g_i)$
corresponding to $h^*$ in $H$.
\begin{figure}[htbp]
  \begin{center}
    \includegraphics[scale=0.7]{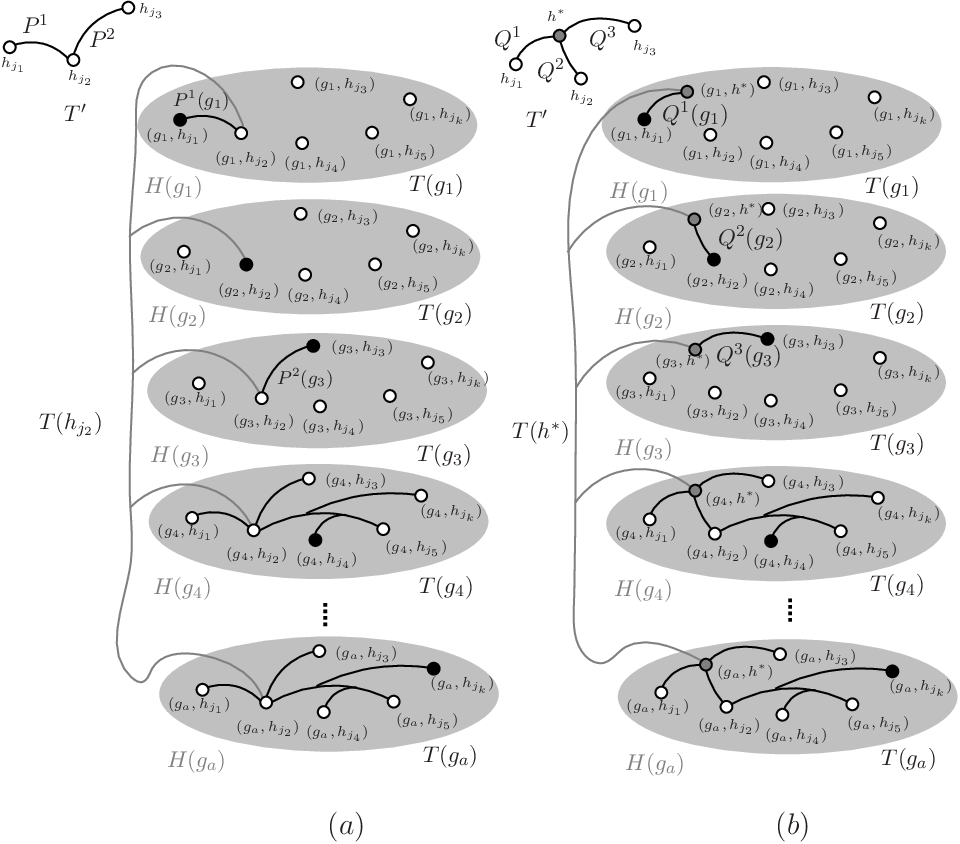}
    \caption{Graphs for Subcase 2.1 in the proof of Theorem
\ref{th2-1}.}
    \label{fig:logo}
  \end{center}
\end{figure}

On the other hand, since there is an $S_G$-Steiner tree of size
$d_G(S_G)$ in $G$, it follows that there exists an Steiner tree of
size $d_G(S_G)$ connecting $\{(g_1,h^*),(g_2,h^*),$
$\ldots,(g_k,h^*)\}$ in $G(h^*)$, say $T(h^*)$. Furthermore, the
subgraph induced by the edges in
$$
\left(\bigcup_{i=4}^{a}E(T(g_i))\right)\cup E(Q^1(g_1))\cup
E(Q^2(g_2))\cup E(Q^3(g_3))\cup E(T(h^*))
$$
is an $S$-Steiner tree in $G\Box H$ (see Fig.3 $(b)$), and hence
$d_{G\Box H}(S)\leq d_G(S_G)+(a-2)d_H(S_H)$. Since $r=a-3$, it
follows that $d_{G\Box H}(S)\leq d_G(S_G)+(r+1)d_H(S_H)$.

{\bf Subcase 2.2.} $x=1$ or $x=2$.

Without loss of generality, let $|V(H(g_{1}))\cap S|=1$ and
$(g_{i_1},h_{j_1})=(g_1,h_{j_1})$. Since there is an $S_G$-Steiner
tree of size $d_G(S_G)$ in $G$, it follows that there exists an
Steiner tree of size $d_G(S_G)$ connecting
$\{(g_1,h_{j_1}),(g_2,h_{j_1}),\ldots,(g_a,h_{j_1})\}$ in
$G(h_{j_1})$, say $T(h_{j_1})$. Since there is an $S_H$-Steiner tree
of size $d_H(S_H)$ in $H$, it follows that there exists an Steiner
tree of size $d_H(S_H)$ connecting
$\{(g_{2},h_{j_1}),(g_2,h_{j_2}),$ $\ldots,(g_2,h_{j_k})\}$ in
$H(g_2)$, say $T(g_2)$. For each $i \ (3\leq i\leq a)$, let $T(g_i)$
be the Steiner tree in $H(g_i)$ corresponding to $T(g_2)$ in
$H(g_2)$. Note that $T(g_i) \ (2\leq i\leq a)$ is the Steiner tree
of size $d_H(S_H)$ connecting
$\{(g_{i},h_{j_1}),(g_i,h_{j_2}),\ldots,(g_i,h_{j_k})\}$ in
$H(g_i)$. Furthermore, the subgraph induced by the edges in
$\left(\bigcup_{i=2}^{a}E(T(g_i))\right)\cup E(T(h_{j_1}))$ is an
$S$-Steiner tree, and hence $d_{G\Box H}(S)\leq
d_G(S_G)+(a-1)d_H(S_H)$. From the definition of $r$, we have $r=a-2$
or $r=a-1$, and hence $d_{G\Box H}(S)\leq d_G(S_G)+(r+1)d_H(S_H)$,
as desired.

From the above argument, we conclude that $d_{G\Box H}(S)\leq
\min\{d_G(S_G)+(r+1)d_H(S_H),$ $d_H(S_H)+(t+1)d_G(S_G)\}$, as
desired. \qed
\end{pf}

\vskip 0.5cm

The following corollaries are immediate from Theorem \ref{th2-1}.
\begin{cor}\label{cor2-1}
Let $G,H$ be two connected graphs of order $n,m$, respectively. Let
$k$ be an integer with $3\leq k\leq mn$. Let
$S=\{(g_{i_1},h_{j_1}),(g_{i_2},h_{j_2}),\ldots,(g_{i_k},h_{j_k})\}$
be a set of distinct vertices of $G\Box H$. Let
$S_G=\{g_{i_1},g_{i_2},\ldots,g_{i_k}\}$ and
$S_H=\{h_{j_1},h_{j_2},\ldots,h_{j_k}\}$. Then
\begin{eqnarray*}
d_G(S_G)+d_H(S_H)&\leq &d_{G\Box H}(S)\\[0.1cm]
&\leq&\min\{d_G(S_G)+(k-2)d_H(S_H),d_H(S_H)+(k-2)d_G(S_G)\}\\[0.1cm]
&=&d_G(S_G)+d_H(S_H)+(k-3)\min\{d_H(S_H),d_G(S_G)\}.
\end{eqnarray*}
\end{cor}

\begin{cor}\label{cor2-2}
Let $G,H$ be two connected graphs, and let $(g,h)$, $(g',h')$ and
$(g'',h'')$ be three vertices of $G\Box H$. Let $S_G=\{g,g',g''\}$,
$S_H=\{h,h',h''\}$, and $S=\{(g,h),(g',h'),$ $(g'',h'')\}$. Then
$$
d_{G\Box H}(S)= d_{G}(S_G)+ d_{H}(S_H)
$$
\end{cor}

To show the sharpness of the above upper and lower bound, we
consider the following example.

\noindent{\bf Example 1:} $(1)$ For $k=3$, from Corollary
\ref{cor2-2}, we have $d_{G\Box H}(S)=d_G(S_G)+d_H(S_H)$, which
implies that the upper and lower bounds in Corollary \ref{cor2-1}
and Theorem \ref{th2-1} are sharp.

$(2)$ Let $G=P_n$ and $H=K_{1,m-1}$, where $P_n=g_1g_2\cdots g_n$,
$h_1,h_2,\cdots, h_{m-1}$ are the leaves of $H$, and $h_{m}$ is the
center of $H$. Choose $S=\{(g_1,h_1),(g_1,h_2),(g_1,h_m)\}\cup
\{(g_n,h_1),(g_n,h_2),(g_n,h_m)\}\cup
\{(g_i,h_1),(g_i,h_2),(g_i,h_m)\,|\,2\leq i\leq x-2\}$, where $4\leq
x\leq n$. Then $d_G(S_G)=n-1$, $d_H(S_H)=2$, $r=x-1$, $t=3$ and
$d_{G\Box H}(S)=n-1+2x=n-1+2+\min\{2(x-1),3(n-1)\}
=d_G(S_G)+d_H(S_H)+\min\{rd_H(S_H),td_G(S_G)\}$, which implies that
the upper bound in Corollary \ref{cor2-1} are sharp.

\subsection{Steiner diameter of Cartesian product graphs}

For Steiner $k$-diameter, we have the following.
\begin{thm}\label{th2-2}
Let $k,m,n$ be an integer with $3\leq k\leq mn$ and $n\leq m$. Let
$G,H$ be two connected graphs of order $n,m$, respectively.

$(1)$ If $k\leq n$, then
\begin{eqnarray*}
&&sdiam_k(G)+sdiam_k(H)\\
&\leq&sdiam_k(G\Box H)\\[0.1cm]
&\leq&sdiam_k(G)+sdiam_k(H)+(k-3)\min\{sdiam_k(G),sdiam_k(H)\}.
\end{eqnarray*}

$(2)$ If $n<k\leq m$, then
\begin{eqnarray*}
n-1+sdiam_k(H)
&\leq&sdiam_k(G\Box H)\\[0.1cm]
&\leq&n-1+sdiam_k(H)+(k-3)\min\{n-1,sdiam_k(H)\}.
\end{eqnarray*}

$(3)$ If $m<k\leq mn$, then
$$
n+m-2\leq sdiam_k(G\Box H)\leq m-1+(k-2)(n-1).
$$

$(4)$ If $mn-\kappa(G\Box H)+1\leq k\leq mn$, then $sdiam_k(G\Box
H)=k-1$.
\end{thm}
\begin{pf}
We first consider all the upper bounds in this theorem. From the
definition of $sdiam_k(G\Box H)$, there exists a vertex subset
$S\subseteq V(G\Box H)$ with $|S|=k$ such that $d_{G\Box
H}(S)=sdiam_k(G\Box H)$. Let
$S=\{(g_{i_1},h_{j_1}),(g_{i_2},h_{j_2}),\ldots,(g_{i_k},h_{j_k})\}$,
and let $S_G=\{g_{i_1},g_{i_2},\ldots,g_{i_k}\}$ and
$S_H=\{h_{j_1},h_{j_2},\ldots,h_{j_k}\}$. From Corollary
\ref{cor2-1}, we have
$$
sdiam_k(G\Box H)=d_{G\Box H}(S)\leq
\min\{d_G(S_G)+(k-2)d_H(S_H),(k-2)d_G(S_G)+d_H(S_H)\}.
$$
For $(1)$, since $k\leq n$, it follows that $d_G(S_G)\leq
sdiam_k(G)$ and $d_H(S_H)\leq sdiam_k(H)$, and hence
\begin{eqnarray*}
&&sdiam_k(G\Box H)=d_{G\Box H}(S)\\
&\leq&\min\{d_G(S_G)+(k-2)d_H(S_H),(k-2)d_G(S_G)+d_H(S_H)\}\\
&\leq&\min\{sdiam_k(G)+(k-2)sdiam_k(H),(k-2)sdiam_k(G)+sdiam_k(H)\}\\
&=&sdiam_k(G)+sdiam_k(H)+(k-3)\min\{sdiam_k(G),sdiam_k(H)\}.
\end{eqnarray*}
For $(2)$, since $n<k\leq m$, it follows that $d_G(S_G)\leq n-1$ and
$d_H(S_H)\leq sdiam_k(H)$, and hence
\begin{eqnarray*}
sdiam_k(G\Box H)&=&d_{G\Box H}(S)\\
&\leq&\min\{d_G(S_G)+(k-2)d_H(S_H),(k-2)d_G(S_G)+d_H(S_H)\}\\
&\leq&\min\{n-1+(k-2)sdiam_k(H),(k-2)(n-1)+sdiam_k(H)\}\\
&=&n-1+sdiam_k(H)+(k-3)\min\{n-1,sdiam_k(H)\}.
\end{eqnarray*}
For $(3)$, since $m<k\leq mn$, it follows that $d_G(S_G)\leq n-1$
and $d_H(S_H)\leq m-1$, and hence
\begin{eqnarray*}
sdiam_k(G\Box H)&=&d_{G\Box H}(S)\\
&\leq&\min\{d_G(S_G)+(k-2)d_H(S_H),(k-2)d_G(S_G)+d_H(S_H)\}\\
&\leq&\min\{n-1+(k-2)(m-1),(k-2)(n-1)+(m-1)\}\\
&=&m-1+(k-2)(n-1).
\end{eqnarray*}

Next, we consider the lower bounds in this theorem. For $(1)$, we
suppose $k\leq n\leq m$. From the definition of $sdiam_k(G)$, it
follows that there exists a vertex subset $S_G\subseteq V(G)$ with
$|S_G|=k$ such that $d_G(S_G)=sdiam_k(G)$. Similarly, there exists a
vertex subset $S_H\subseteq V(H)$ with $|S_H|=k$ such that
$d_H(S_H)=sdiam_k(H)$. Without loss of generality, let
$S_G=\{g_1,g_2,\ldots,g_k\}$ and $S_H=\{h_1,h_2,\ldots,h_k\}$. Then
$S=\{(g_1,h_1),(g_2,h_2),\ldots,(g_k,h_k)\} \subseteq V(G\Box H)$
and $|S|=k$. From Lemma \ref{lem2-2} and the definition of Steiner
$k$-diameter, we have
$$
sdiam_k(G)+sdiam_k(H)=d_G(S_G)+d_H(S_H)\leq d_{G\Box H}(S)\leq
sdiam_k(G\Box H).
$$

For $(2)$, we suppose $n<k\leq m$. Let
$S=\{(g_{i_1},h_{j_1}),(g_{i_2},h_{j_2}),\ldots,(g_{i_k},h_{j_k})\}$
be a set of distinct vertices of $G\Box H$ such that $V(G)\subseteq
\{g_{i_1},g_{i_2},\ldots,g_{i_k}\}=S_G$ and $d_H(S_H)=sdiam_k(H)$,
where $S_H=\{h_{j_1},h_{j_2},\ldots,h_{j_k}\}$. From Lemma
\ref{lem2-2}, we have
$$
n-1+sdiam_k(H)=d_G(S_G)+d_H(S_H)\leq d_{G\Box H}(S)\leq
sdiam_k(G\Box H).
$$

For $(3)$, we suppose $m<k\leq mn$. Let
$S=\{(g_{i_1},h_{j_1}),(g_{i_2},h_{j_2}),\ldots,(g_{i_k},h_{j_k})\}$
be a set of distinct vertices of $G\Box H$ such that $V(G)\subseteq
S_G$ and $V(H)\subseteq S_H$, where
$S_G=\{g_{i_1},g_{i_2},\ldots,g_{i_k}\}$ and
$S_H=\{h_{j_1},h_{j_2},\ldots,h_{j_k}\}$. From Lemma \ref{lem2-2},
we have
$$
n+m-2=(n-1)+(m-1)=d_G(S_G)+d_H(S_H)\leq d_{G\Box H}(S)\leq
sdiam_k(G\Box H),
$$
as desired.

For $(4)$, we suppose $mn-\kappa(G\Box H)+1\leq k\leq mn$. For any
$S\subseteq V(G\Box H)$ with $|S|=k$, we have $|V(G)|-|S|\leq
\kappa(G\Box H)-1$, and hence $G[S]$ is connected. Therefore, we
have $d_{G\Box H}(S)\leq k-1$, and hence $sdiam_k(G\Box H)\leq k-1$
by the arbitrariness of $S$. So, we have $sdiam_k(G\Box H)=k-1$.
\qed
\end{pf}

\vskip0.5cm

The following corollary is immediate from Theorem \ref{th2-2}.
\begin{cor}\label{cor2-3}
Let $G,H$ be two connected graphs of order at least $3$. Then
$$
sdiam_3(G\Box H)=sdiam_3(G)+sdiam_3(H).
$$
\end{cor}

To show the sharpness of the above upper and lower bound, we
consider the following example.

\noindent{\bf Example 2:} $(1)$ For $k=3$, from Corollary
\ref{cor2-3}, we have $sdiam_k(G\Box H)=sdiam_k(G)+sdiam_k(H)$,
which implies that the upper and lower bounds in Theorem \ref{th2-2}
are sharp.

$(2)$ Let $G=P_n$ and $H=P_m$ with $5\leq n\leq m$. Then
$sdiam_4(G)=n-1$, $sdiam_4(H)=m-1$ and $sdiam_4(G\Box
H)=2(n-1)+(m-1)$, which implies that all the upper bounds in Theorem
\ref{th2-2} are sharp.

\section{Results for lexicographic product}

From the definition, the lexicographic product graph $G\circ H$ is
the graph obtained by replacing each vertex of $G$ by a copy of $H$
and replacing each edge of $G$ by a complete bipartite graph
$K_{m,m}$, where $m=|V(H)|$.
\begin{lem}{\upshape \cite{Hammack}}\label{lem3-1}
Let $G$ and $H$ be two graphs, and let $(g,h)$ and $(g',h')$ be two
vertices of $G\circ H$. Then
$$
d_{G\circ H}((g,h),(g',h'))=\left\{
\begin{array}{ll}
d_{G}(g,g'), &\mbox {\rm if} \ g\neq g';\\[0.2cm]
d_{H}(h,h'), &\mbox {\rm if} \ g=g'~and~deg_{G}(g)=0;\\[0.2cm]
\min \{d_{H}(h,h'),2\}, &\mbox {\rm if} \ g=g'~and~deg_{G}(g)\neq 0.
\end{array}
\right.
$$
\end{lem}

A \emph{weak homomorphism} $\varphi: G\rightarrow H$ is a map
$\varphi: V(G)\rightarrow V(H)$ for which $uv\in E(G)$ implies
$\varphi(u)\varphi(v)\in E(H)$ or $\varphi(u)=\varphi(v)$. Observe
that the projection $p:G\circ H\rightarrow G$ is a weak
homomorphism. For more details, we refer to \cite{Hammack}
(p.32,p.57).
\begin{lem}{\upshape \cite{Hammack}}\label{lem3-2}
Let $G$ and $H$ be two graphs, and let $(g,h)$ and $(g',h')$ be two
vertices of $G\circ H$. Then
$$
d_{G\circ H}((g,h),(g',h'))\geq d_{G}(g,g').
$$
\end{lem}

\subsection{Steiner distance of lexicographic product graphs}

The following lemma is a generalization of Lemma \ref{lem3-2}, which
is a natural lower bound of $d_{G\circ H}(S)$ for $S\subseteq
V(G\circ H)$ and $|S|=k$.
\begin{lem}\label{lem3-3}
Let $k\geq 2$ be an integer, $G$ be a connected graph, and $H$ be a
graph. Let $S=\{(g_{i_1},h_{j_1}),$
$(g_{i_2},h_{j_2}),\ldots,(g_{i_k},h_{j_k})\}$ be a set of distinct
vertices of $G\circ H$. Let
$S_G=\{g_{i_1},g_{i_2},\ldots,g_{i_k}\}$. Then
$$
d_{G\circ H}(S)\geq d_{G}(S_G).
$$
\end{lem}
\begin{pf}
We note that $g_{i_1},g_{i_2},\ldots,g_{i_k}$ are not necessarily
distinct. Let $T$ be a minimum $S$-Steiner tree in $G\circ H$. So
$T$ has $d_{G\circ H}(S)$ edges. Let $Z$ be the minor obtained from
$G\circ H$ by contracting edge in $H(g)$ for every $g$ of $G$.
(Equivalently, identifying all the vertices in $H(g)$ into a single
vertex $g$ and delete multiple edges in the resulting graph.) Then
$Z$ is isomorphic to $G$. Now $T$ becomes $Y$, a connected subgraph
of $Z$ containing the vertices corresponding to
$g_{i_1},g_{i_2},\ldots,g_{i_k}$ in $G$. Thus $E(Y)\geq d_G(S_G)$.
Since $E(T)\geq E(Y)$, the result follows.\qed
\end{pf}

\cite{AnandCKP} obtained the following formula.
\begin{lem}{\upshape \cite{AnandCKP}}\label{lem3-4}
Let $k\geq 2$. Let $G,H$ be two graphs such that $G$ is connected.
Let $S=\{(g_{i_1},h_{j_1}),$
$(g_{i_2},h_{j_2}),\ldots,(g_{i_k},h_{j_k})\}$ be a set of distinct
vertices of $G\circ H$ such that $g_{i_p}\neq g_{i_q} \ (1\leq
p,q\leq k)$. Let $S_G=\{g_{i_1},g_{i_2},\ldots,g_{i_k}\}$. Then
$$
d_{G\circ H}(S)=d_{G}(S_G).
$$
\end{lem}

For general case, we have the following formula for Steiner distance
of lexicographic product graphs.
\begin{thm}\label{th3-1}
Let $k,n,m$ be three integers with $2\leq k\leq mn$. Let $G$ be a
connected graph of order $n$, and $H$ be a graph of order $m$. Let
$S=\{(g_{i_1},h_{i_1}),(g_{i_2},h_{i_2}),\ldots,(g_{i_k},h_{i_k})\}$
be a set of distinct vertices of $G\circ H$. Let
$S_G=\{g_{i_1},g_{i_2},\ldots,g_{i_k}\}$ and
$S_H=\{h_{j_1},h_{j_2},\ldots,h_{j_k}\}$ (note that $S_G,S_H$ are
both multi-sets). Let $r$ be the number of distinct vertices in
$S_G$, where $1\leq r\leq k$.

$(1)$ If $r=1$ and $H[S_H]$ is connected in $H$, then $d_{G\circ
H}(S)=k-1$.

$(2)$ If $r=1$ and $H[S_H]$ is not connected in $H$, then $d_{G\circ
H}(S)=k$.

$(3)$ If $r\geq 2$, then $d_{G\circ H}(S)=d_{G}(S_G)+k-r$.
\end{thm}
\begin{pf}
$(1)$ Since $r=1$, it follows that $g_{i_1}=g_{i_2}=\ldots=g_{i_k}$,
and hence $S=\{(g_{i_1},h_{i_1}),(g_{i_2},h_{i_2}),\ldots,$
$(g_{i_k},h_{i_k})\}=\{(g_{i_1},h_{i_1}),(g_{i_1},h_{i_2}),\ldots,(g_{i_1},h_{i_k})\}\subseteq
V(H(g_{i_1}))$. Since $H[S_H]$ is connected in $H$, it follows that
the subgraph induced by the vertices in $\{(g_{i_1},h_{i_1}),$
$(g_{i_2},h_{i_2}),\ldots,(g_{i_k},h_{i_k})\}$ is connected in
$H(g_{i_1})$, and hence $d_{G\circ H}(S)=k-1$.

$(2)$ Since $H[S_H]$ is not connected in $H$, it follows that the
subgraph induced by the vertices in $\{(g_{i_1},h_{i_1}),$
$(g_{i_2},h_{i_2}),\ldots,(g_{i_k},h_{i_k})\}$ is not connected in
$H(g_{i_1})$, and hence $d_{G\circ H}(S)\geq k$. Since $G$ is a
connected graph of order at least $2$, it follows that there exists
a vertex $g^*\in V(G)$ such that $g_{i_1}g^*\in E(G)$. From the
structure of $G\circ H$, the tree induced by the edges in
$\{(g_{i_p},h_{i_p})(g^*,h_1)\,|\,1\leq p\leq
k\}=\{(g_{i_1},h_{i_p})(g^*,h_1)\,|\,1\leq p\leq k\}$ is an
$S$-Steiner tree in $G\circ H$, and hence $d_{G\circ H}(S)\leq k$.
So, we have $d_{G\circ H}(S)=k$.

$(3)$ Since $r\geq 2$, it follows that the vertices in $S$ belong to
at least two copies of $H$ in $G\circ H$. From the definition of
$r$, we can assume that $H(g_1),H(g_2),\ldots,H(g_{r})$ satisfy
$S\cap V(H(g_i))\neq \emptyset$ for each $g_i \ (1\leq i\leq r)$,
and $S\cap V(H(g_i))=\emptyset$ for each $g_i \ (r+1\leq i\leq n)$.
Let $S_G'=\{g_1,g_2,\ldots,g_{r}\}$. Then $S_G'=S_G$ when we regard
$S_G$ as a normal set, not a multi-set. Clearly,
$d_G(S_G)=d_G(S_G')$, and
$S=\{(g_{i_1},h_{i_1}),(g_{i_2},h_{i_2}),\ldots,(g_{i_k},h_{i_k})\}\subseteq
\bigcup_{i=1}^rV(H(g_i))$. Without loss of generality, we can assume
$(g_{i_a},h_{i_a})\in V(H(g_a))$ for each $a \ (1\leq a\leq r)$.
Then $(g_{i_a},h_{i_a})=(g_{a},h_{i_a})$ for each $a \ (1\leq a\leq
r)$. Let $S'=\{(g_{a},h_{i_a})\,|\,1\leq a\leq r\}$. Then
$(g_{i_{r+1}},h_{i_{r+1}}),(g_{i_{r+2}},h_{i_{r+2}}),\ldots,(g_{i_k},h_{i_k})\in
(\bigcup_{i=1}^rV(H(g_i))\setminus S'$. Note that there exists an
$S_G'$-Steiner tree $T_G$ of size $d_{G}(S_G')=d_{G}(S_G)$ in $G$.
Without loss of generality, let $V(T_G)=\{g_1,g_2,\ldots,g_{t}\}$,
where $r\leq t\leq n$. In order to select $d_{G}(S_G')$ edges in
$G\circ H$ to form an $S'$-Steiner tree $T'$ in $G\circ H$
isomorphic to $T_G$ in $G$ such that $V(T')\subseteq
\bigcup_{i=1}^tV(H(g_i))$, we define a function $f:
E(T_G)\longrightarrow E(T')$:
$$
f(g_ag_b)=\left\{
\begin{array}{ll}
(g_a,h_{i_a})(g_b,h_{i_b}), &\mbox {\rm if}~g_a\in S_G~{\rm and}~g_b\in S_G;\\[0.2cm]
(g_a,h_{i_a})(g_b,h_{1}), &\mbox {\rm if}~g_a\in S_G~{\rm and}~g_b\notin S_G;\\[0.2cm]
(g_a,h_{1})(g_b,h_{1}), &\mbox {\rm if}~g_a\notin S_G~{\rm
and}~g_b\notin S_G,
\end{array}
\right.
$$
for each $g_ag_b\in E(T_G) \ (1\leq a\neq b\leq t)$. Note that $T'$
is an $S'$-Steiner tree in $G\circ H$.

We now extend the tree $T'$ to an $S$-Steiner tree $T$ by adding
$|S|-|S'|=k-r$ edges in $G\circ H$. For each vertex
$(g_{i_a},h_{j_a})\in S\setminus S' \ (r+1\leq a\leq k)$, since
there exists a vertex $g_{i_b}\in V(T_G) \ (1\leq b\neq a\leq t)$ in
$G$ such that $g_{i_a}g_{i_b}\in E(T_G)$, we select an edge
$$
e_a=\left\{
\begin{array}{ll}
(g_a,h_{i_a})(g_b,h_{i_b}), &\mbox {\rm if}~g_b\in S_G;\\[0.2cm]
(g_a,h_{i_a})(g_b,h_{1}), &\mbox {\rm if}~g_b\notin S_G
\end{array}
\right.
$$
in $G\circ H$, and then add it into $T'$. Observe that the tree
induced by the edges in $E(T')\cup \{e_{a}\,|\,r+1\leq i\leq k\}$ is
an $S$-Steiner tree $T$ in $G\circ H$. Since $|E(T)|=d_G(S_G)+k-r$,
it follows that $d_{G\circ H}(S)\leq d_{G}(S_G)+k-r$.

It remains us to show that $d_{G\circ H}(S)\geq d_{G}(S_G)+k-r$.
Recall that $V(G)=\{g_1,g_2,\ldots,g_{n}\}$. Without loss of
generality, we assume that $H(g_{1}),H(g_{2}),\ldots,H(g_{r})$ be
the $H$ copies such that $V(H(g_{i}))\cap S\neq \emptyset$, $1\leq
i\leq r$. Clearly,
$S=\{(g_{i_1},h_{i_1}),(g_{i_2},h_{i_2}),\ldots,(g_{i_k},h_{i_k})\}\subseteq
\bigcup_{i=1}^rV(H(g_i))$. Set $|S\cap V(H(g_i))|=x_i$. Then
$\sum_{i=1}^{r}x_i=k$. Without loss of generality, let $S_i=S\cap
V(H(g_i))=\{(g_i,h_j)\,|\,1\leq j\leq x_i\}$ for each $g_i \ (1\leq
i\leq r)$. In order to find an $S$-Steiner tree $T$ in $G\circ H$,
we need the edges between some $H(g_i)$ and $H(g_j)$, $1\leq i\neq
j\leq r$. Note that $S_i\subseteq V(H(g_i))$ for each $i \ (1\leq
i\leq r)$. Clearly, there exists a subtree $T'$ connecting $S'$ in
$T$ such that $E(T')\cap (\bigcup_{i=1}^rE(H(g_i)))=\emptyset$,
where $|S'\cap S_i|=1 \ (1\leq i\leq r)$. Since $|E(T')|\geq
d_G(S_G)$ and $|S|-|S'|=k-r$, it follows that $T$ is an $S$-Steiner
tree of size $d_G(S_G)+k-r$ in $G\circ H$, and hence $d_{G\circ
H}(S)\geq d_{G}(S_G)+k-r$.

From the above argument, we conclude that $d_{G\circ H}(S)=
d_{G}(S_G)+k-r$.\qed
\end{pf}

\vskip 0.3cm

In Theorem \ref{th3-1}, we assume that $G$ is a connected graph. For
$k=3$, we have the following by assuming that $G$ is not connected.
\begin{pro}\label{pro3-1}
Let $G$ and $H$ be two graphs such that $G$ is connected, and let
$(g,h)$, $(g',h')$ and $(g'',h'')$ be three vertices of $G\circ H$.
Let $S=\{(g,h),(g',h'),(g'',h'')\}$, $S_G=\{g,g',g''\}$ and
$S_H=\{h,h',h''\}$. Then
$$
d_{G\circ H}(S)=\left\{
\begin{array}{ll}
d_{H}(S_H), &\mbox {\rm if} \ g=g'=g''~and~deg_{G}(g)=0;\\[0.2cm]
\min \{d_{H}(S_H),3\}, &\mbox {\rm if} \ g=g'=g''~and~deg_{G}(g)\neq 0;\\[0.2cm]
\infty, &\mbox {\rm if} \ g\neq g',~g'=g''~and~d_G(g,g')=\infty;\\[0.2cm]
d_{G}(g,g')+1, &\mbox {\rm if} \ g\neq g',~g'=g''~and~d_G(g,g')\neq \infty;\\[0.2cm]
d_{G}(S_G), &\mbox {\rm if} \ g\neq g',~g\neq g''~and~g'\neq g''.
\end{array}
\right.
$$
\end{pro}
\begin{pf}
Suppose that $g=g'=g''$ and $deg_{G}(g)=0$. Since $g$ is isolated,
it follows that $H(g)$ is a component of $G\circ H$, and hence
$d_{G\circ H}(S)=d_H(S_H)$.

Suppose that $g=g'=g''$ and $deg_{G}(g)\geq 1$. Since
$deg_{G}(g)\geq 1$, there exists a vertex $g^*$ in $G$ such that
$gg^*\in E(G)$, and hence the tree induced by the edges in
$$
\{(g,h)(g^*,h),(g,h')(g^*,h),(g,h'')(g^*,h)\}
$$
is an $S$-Steiner tree. Therefore, $d_{G\circ H}(S)\leq 3$. On the
other hand, from Observation \ref{obs1-1}, $d_{G\circ H}(S)\geq 2$.
So $d_{G\circ H}(S)=2$ or $d_{G\circ H}(S)=3$. Since $d_H(S_H)\geq
2$ by Observation \ref{obs1-1}, it follows that $d_{G\circ H}(S)=
\min\{d_H(S_H),3\}$.

Suppose that $g\neq g'$, $g'=g''$ and $d_G(g,g')=\infty$. Then there
is no path connecting $g$ and $g'$ in $G$. Note that $(g,h)\in
V(H(g))$ and $(g',h'),(g'',h'')\in V(H(g'))$. Clearly, there is no
$S$-Steiner tree in $G\circ H$. Therefore, $d_{G\circ H}(S)=\infty$.

Suppose that $g\neq g'$, $g'=g''$ and $d_G(g,g')\neq \infty$. Set
$d_G(g,g')=\ell$. Let $P=gg_1g_2\cdots g_{\ell-1}g'$ be a path
connecting $g$ and $g'$ in $G$. Then the tree induced by the edges
in
$$
\{(g,h)(g_1,h),(g_1,h)(g_2,h),\cdots,(g_{\ell-2},h)(g_{\ell-1},h),(g_{\ell-1},h)(g',h'),(g_{\ell-1},h)(g'',h'')\}
$$
is an $S$-Steiner tree. Therefore, $d_{G\circ H}(S)\leq \ell+1$. It
suffices to show $d_{G\circ H}(S)\geq \ell+1$. From Observation
\ref{obs2-1}, any minimal $S$-Steiner tree $T$ is a path or there
exists a vertex $(g^*,h^*)\in V(G\circ H)\setminus S$ such that the
tree $T$ consists of three paths connecting $(g^*,h^*)$ and
$(g,h),(g',h'),(g'',h'')$, respectively. If $T$ is a path, then we
can assume that $(g',h')$ be the internal vertex of the path $T$.
Since $g'=g''$, it follows that $(g',h'),(g'',h'')\in V(H(g'))$. One
can see that the length of the path from $(g',h')$ to $(g'',h'')$ is
at least $1$. By Lemma \ref{lem3-2}, $d_{G\circ H}(S)=d_{G\circ
H}((g,h)(g',h'))+1\geq d_{G}(g,g')+1=\ell+1$, as desired. Suppose
that $T$ is a tree and there exists a vertex $(g^*,h^*)\in V(G\circ
H)\setminus S$ such that $T$ consists of three paths connecting
$(g^*,h^*)$ and $(g,h),(g',h'),(g'',h'')$, respectively. Then
\begin{eqnarray*}
d_{G\circ H}(S)&=&d_{G\circ H}((g,h)(g^*,h^*))+d_{G\circ H}((g',h')(g^*,h^*))+d_{G\circ H}((g'',h'')(g^*,h^*))\\
&\geq&d_{G}(g,g^*)+d_{G}(g',g^*)+1 \ \mbox{\rm (by~Lemma \ref{lem3-2})}\\
&\geq&\ell+1.
\end{eqnarray*}

Suppose that $g\neq g'$, $g\neq g''$ and $g'\neq g''$. From Lemma
\ref{lem3-4}, we have $d_{G\circ H}(S)=d_{G}(S_G)$, as desired. The
proof is now complete. \qed
\end{pf}

\subsection{Steiner diameter of lexicographic product graphs}

By Theorem \ref{th3-1}, we can derive the following results for
Steiner diameter of lexicographic product graphs.
\begin{thm}\label{th3-2}
Let $k,n,m$ be three integers with $2\leq k\leq mn$. Let $G$ be a
connected graph of order $n$, and $H$ be a graph of order $m$. Then

$(1)$
$$
sdiam_k(G\circ H)\leq \left\{
\begin{array}{ll}
sdiam_k(G)+k-2, &\mbox {\rm if} \ 2\leq k\leq n;\\[0.2cm]
\max\{n+k-3,k\}, &\mbox {\rm if} \ n<k\leq mn;
\end{array}
\right.
$$
Furthermore, if $n\geq 3$, then $sdiam_k(G\circ H)\leq n+k-3$.

$(2)$
$$
sdiam_k(G\circ H)\geq \left\{
\begin{array}{ll}
sdiam_k(G), &\mbox {\rm if} \ 2\leq k\leq n;\\[0.2cm]
k-1, &\mbox {\rm if} \ n<k\leq nm.
\end{array}
\right.
$$
Moreover, if
$$
r=\min_{2\leq x\leq n}\{x\,|\,sdiam_x(G)=n-1\},
$$
then
\begin{eqnarray*}
&&sdiam_k(G\circ H)\geq\\[0.2cm]
&&\left\{
\begin{array}{ll}
sdiam_k(G), &\mbox {\rm if} \ 2\leq k\leq r;\\[0.2cm]
n-1+k-r, &\mbox {\rm if} \ r<k\leq rm;\\[0.2cm]
n-1+r(m-1)+\lfloor \frac{k-rm}{m}\rfloor(m-1)+\max\{k-(r+\lfloor
\frac{k-rm}{m}\rfloor)m-1,0\}, &\mbox {\rm if} \ rm<k\leq nm.
\end{array}
\right.
\end{eqnarray*}
\end{thm}
\begin{pf}
$(1)$ From the definition of $sdiam_k(G\circ H)$, there exists a
vertex subset $S\subseteq V(G\circ H)$ with $|S|=k$ such that
$d_{G\circ H}(S)=sdiam_k(G\circ H)$. Let
$S=\{(g_{i_1},h_{i_1}),(g_{i_2},h_{i_2}),\ldots,(g_{i_k},h_{i_k})\}$,
and $S_G=\{g_{i_1},g_{i_2},\ldots,g_{i_k}\}$. Let $s$ be the number
of distinct vertices in $S_G$. We apply Theorem \ref{th3-1}. (Here
$s$ plays the role of $r$ in Theorem \ref{th3-1}.) If $s\geq 2$,
then $sdiam_k(G\circ H)=d_{G\circ H}(S)=d_{G}(S_G)+k-s\leq
d_{G}(S_G)+k-2$. Furthermore, if $k\leq n$, then $d_{G}(S_G)\leq
sdiam_k(G)$, and hence $sdiam_k(G\circ H)\leq d_{G}(S_G)+k-2\leq
sdiam_k(G)+k-2$. If $n<k\leq mn$, then $d_{G}(S_G)\leq n-1$, and
hence $sdiam_k(G\circ H)\leq d_{G}(S_G)+k-2\leq (n-1)+k-2=n+k-3$.
Note that if $s=1$, then $k\leq m$, and hence $sdiam_k(G\circ
H)=d_{G\circ H}(S)\leq k$. From the above argument, we conclude that
$sdiam_k(G\circ H)\leq sdiam_k(G)+k-2$ if $k\leq n$, and
$sdiam_k(G\circ H)\leq \max\{n+k-3,k\}$ if $n<k\leq mn$, as desired.

$(2)$ If $2\leq k\leq n$, then we let $S=\{(g_{i_1},h_{i_1}),$
$(g_{i_2},h_{i_2}),\ldots,(g_{i_k},h_{i_k})\}$ be a set of distinct
vertices of $G\circ H$ such that $d_G(S_G)=sdiam_k(G)$, where
$S_G=\{g_{i_1},g_{i_2},\ldots,g_{i_k}\}$. From Lemma \ref{lem3-3},
we have $sdiam_k(G)=d_G(S_G)\leq d_{G\circ H}(S)\leq sdiam_k(G\circ
H)$. If $n\leq k\leq nm$, then it follows from Observation
\ref{obs1-1} that $k-1\leq d_{G\circ H}(S)\leq sdiam_k(G\circ H)$
for any $S\subseteq V(G\circ H)$ and $|S|=k$.

Now for the ``moreover'' part of the result. Let $r=\min_{2\leq
x\leq n}\{x\,|\,sdiam_x(G)=n-1\}$. Suppose $sdiam_r(G)=n-1 \ (2\leq
r\leq n)$. If $2\leq k\leq r$, then $2\leq k\leq r\leq n$, and hence
$sdiam_k(G\circ H)\geq sdiam_k(G)$. Suppose $r<k\leq rm$. Since
$sdiam_r(G)=n-1$, it follows that there exists a vertex set
$S'=\{g_1,g_2,\ldots,g_r\}\subseteq V(G)$ such that
$d_G(S')=n-1=sdiam_k(G)$. Let $S=S_1\cup S_2\subseteq
\bigcup_{i=1}^rV(H(g_i))$ such that $S_1=\{(g_i,h_1)\,|\,1\leq i\leq
r\}$ and $S_2\subseteq \bigcup_{i=1}^rV(H(g_i))-S_1$ and
$|S_2|=k-r$. Since $r\geq 2$ and $sdiam_r(G)=n-1$, it follows that
$sdiam_k(G\circ H)\geq d_{G\circ
H}(S)=d_G(S_G)+k-r=d_G(S')+k-r=n-1+k-r$, as desired. Suppose
$rm<k\leq nm$. Since $sdiam_r(G)=n-1$, it follows that there exists
a vertex set $S'=\{g_1,g_2,\ldots,g_r\}\subseteq V(G)$ such that
$d_G(S')=n-1=sdiam_k(G)$. Let $S=S_1\cup S_2\subseteq
\bigcup_{i=1}^rV(H(g_i))$ such that $S_1=\bigcup_{i=1}^rV(H(g_i))$
and $S_2\subseteq \bigcup_{i=r+1}^xV(H(g_i))$ and $|S_2|=k-rm$,
where $x=\lceil \frac{k-rm}{m}\rceil$. Then $sdiam_k(G\circ H)\geq
d_{G\circ H}(S)=d_G(S')+r(m-1)+\lfloor
\frac{k-rm}{m}\rfloor(m-1)+\max\{k-(r+\lfloor
\frac{k-rm}{m}\rfloor)m-1,0\}=n-1+r(m-1)+\lfloor
\frac{k-rm}{m}\rfloor(m-1)+\max\{k-(r+\lfloor
\frac{k-rm}{m}\rfloor)m-1,0\}$. \qed
\end{pf}

\vskip0.5cm

To show the sharpness of the upper and lower bounds in Theorem
\ref{th3-2}, we consider the following example.

\vskip0.3cm

\noindent{\bf Example 3:} Let $G=P_n$, and $H$ be a graph of order
$m$. If $k\leq \min\{2m,n\}$, then $sdiam_k(G\circ
H)=n+k-3=sdiam_k(G)+k-2$. If $\max\{n,m+1\}\leq k\leq 2m$, then
$sdiam_k(G\circ H)=n+k-3=\max\{n+k-3,k\}$. These implies that the
upper bounds in Theorem \ref{th3-2} are sharp.

\vskip0.3cm

\noindent{\bf Example 4:} Let $G=K_n$ and $H=K_m$. Then $G\circ H$
is a complete graph of order $mn$. If $2\leq k\leq n$, then
$sdiam_k(G)=k-1=sdiam_k(G\circ H)$. If $n\leq k\leq nm$, then
$sdiam_k(G\circ H)=k-1$. These implies that the lower bounds in
Theorem \ref{th3-2} are sharp.

\vskip0.3cm

\noindent{\bf Example 5:} Let $G=P_n \ (n\geq 3)$, and $H$ be a
graph of order $m$. From the definition of $r$, we have $r=2$. For
$2\leq k\leq r$, we have $k=r=2$, and hence
$sdiam_2(G)=n-1=sdiam_2(G\circ H)$. For $r<k\leq rm$, we have
$n-1+k-2\leq sdiam_k(G\circ H)\leq n+k-3$, and hence $sdiam_k(G\circ
H)=n+k-3$. Let $G'=P_n \ (n\geq 3)$, and $H'=P_2$. For $rm<k\leq
nm$, we let $k=2t$. From Theorem \ref{th3-2}, we have
$sdiam_k(G\circ H)\geq n-1+t$. One can easily check that
$sdiam_k(G\circ H)=n-1+t$. These implies that the lower bounds for
parameter $r$ in Theorem \ref{th3-2} are sharp.

The following result is immediate from Proposition \ref{pro3-1}.
\begin{pro}\label{pro3-5}
Let $G,H$ be two connected graphs. Then
$$
sdiam_3(G\circ H)=\left\{
\begin{array}{ll}
diam(G)+1 &\mbox {\rm if} \ G=P_n, \ diam(G)\geq 2,\\[0.2cm]
sdiam_3(G) &\mbox {\rm if} \ G\neq P_n, \ diam(G)\geq 2,\\[0.2cm]
\min \{sdiam_3(H),3\} &\mbox {\rm if} \ G=K_n.
\end{array}
\right.
$$
\end{pro}

\section{Applications}

In this section, we demonstrate the usefulness of the proposed
constructions by applying them to some instances of Cartesian and
lexicographical product networks.

The following results are immediate.
\begin{pro}\label{pro4-1}
Let $k,n$ be two integers with $2\leq k\leq n$.

$(1)$ For a complete graph $K_n$, $sdiam_k(K_n)=k-1$;

$(2)$ For a path $P_n$, $sdiam_k(P_n)=n-1$;

$(3)$ For a cycle $C_n$, $sdiam_k(C_n)=\big
\lfloor\frac{n(k-1)}{k}\big \rfloor$.
\end{pro}

\subsection{Two-dimensional grid graph}

A \emph{two-dimensional grid graph} $G_{n,m}$ is the Cartesian
product graph $P_n\Box P_m$ of path graphs on $m$ and $n$ vertices.
For more details on grid graph, we refer to \cite{Calkin, Itai}. The
network $P_n\circ P_m$ is the lexicographical product of $P_n$ and
$P_m$; see \cite{Mao2}.
\begin{pro}\label{pro4-2}
Let $k,n,m$ be three integers with $3\leq k\leq mn$, $n\geq 3$, and
$m\geq 3$.

$(1)$ For network $P_n\Box P_m$,
$$
m+n-2\leq sdiam_k(P_n\Box P_m)\leq m+n-2+(k-3)\min\{m-1,n-1\}.
$$

$(2)$ For network $P_n\circ P_m$,
$$
n+k-3\geq sdiam_k(P_n\circ P_m)\geq \left\{
\begin{array}{ll}
k-1, &\mbox {\rm if} \ n+1\leq k\leq mn;\\[0.2cm]
n-1, &\mbox {\rm if} \ 2\leq k\leq n.
\end{array}
\right.
$$
\end{pro}
\begin{pf}
$(1)$ From $(2)$ of Proposition \ref{pro4-1}, we have
$sdiam_k(P_n)=n-1$ and $sdiam_k(P_m)=m-1$. By Theorem \ref{th2-2},
$sdiam_k(P_n\Box P_m)\geq sdiam_k(P_n)+sdiam_k(P_m)=m+n-2$ and
$sdiam_k(P_n\Box P_m)\leq m+n-2+(k-3)\min\{m-1,n-1\}$.

$(2)$ Set $G=P_n$ and $H=P_m$. From Theorem \ref{th3-2}, the result
holds.\qed
\end{pf}

\subsection{$r$-dimensional mesh}

An \emph{$r$-dimensional mesh} is the Cartesian product of $r$
paths. By this definition, two-dimensional grid graph is a
$2$-dimensional mesh. An $r$-dimensional hypercube is a special case
of an $r$-dimensional mesh, in which the $r$ linear arrays are all
of size $2$; see \cite{Johnsson}.
\begin{pro}\label{pro4-3}
Let $k,m_1,m_2,\cdots,m_r$ be the integers with $m_1\geq m_2 \geq
\cdots \geq m_r$ and $3\leq k\leq \prod_{i=1}^rm_i$.

$(1)$ For an $r$-dimensional mesh $P_{m_1}\Box P_{m_2}\Box \cdots
\Box P_{m_r}$,
$$
\sum_{i=1}^{r}m_i-r\leq sdiam_k(P_{m_1}\Box P_{m_2}\Box \cdots \Box
P_{m_r})\leq (k-2)\left(\sum_{i=2}^{r}m_i-r+1\right)+m_1-1.
$$

$(2)$ For an $r$-dimensional network $P_{m_1}\circ P_{m_2}\circ
\cdots \circ P_{m_r}$,
$$
m_1+k-2\geq sdiam_k(P_{m_1}\circ P_{m_2}\circ \cdots \circ
P_{m_r})\geq \left\{
\begin{array}{ll}
k-1, &\mbox {\rm if} \ m_1+1\leq k\leq \prod_{i=1}^rm_i;\\[0.2cm]
m_1-1, &\mbox {\rm if} \ 2\leq k\leq m_1.
\end{array}
\right.
$$
\end{pro}
\begin{pf}
$(1)$ From $(2)$ of Proposition \ref{pro4-1},
$sdiam_k(P_{m_i})=m_i-1$ for each $i \ (1\leq i\leq r)$. From
Theorem \ref{th2-2}, we have $sdiam_k(P_{m_1}\Box P_{m_2}\Box \cdots
\Box P_{m_r})\geq
\sum_{i=1}^{r}sdiam_k(P_{m_i})=\sum_{i=1}^{r}m_i-r$, and
$sdiam_k(G\Box H)\leq
sdiam_k(G)+sdiam_k(H)+(k-3)\min\{sdiam_k(G),sdiam_k(H)\}$ for two
connected graphs $G$ and $H$, and hence
\begin{eqnarray*} \label{n-t=0}
&&sdiam_k(P_{m_1}\Box P_{m_2}\Box \cdots \Box P_{m_r})\\
&=&sdiam_k((P_{m_1}\Box P_{m_2}\Box \cdots \Box P_{m_{r-1}})\Box P_{m_r})\\
&\leq&sdiam_k(P_{m_1}\Box P_{m_2}\Box \cdots \Box P_{m_{r-1}})+sdiam_k(P_{m_r})\\
&&+(k-3)\min\{sdiam_k(P_{m_1}\Box P_{m_2}\Box \cdots \Box P_{m_{r-1}}),sdiam_k(P_{m_r})\}\\
&=&sdiam_k(P_{m_1}\Box P_{m_2}\Box \cdots \Box P_{m_{r-1}})+(k-2)(m_{r}-1)\\
&\leq&sdiam_k(P_{m_1}\Box P_{m_2}\Box \cdots \Box P_{m_{r-2}})+sdiam_k(P_{m_{r-1}})+(k-2)(m_{r}-1)\\
&&+(k-3)\min\{sdiam_k(P_{m_1}\Box P_{m_2}\Box \cdots \Box P_{m_{r-2}}),sdiam_k(P_{m_{r-1}})\}\\
\end{eqnarray*}
\begin{eqnarray*}
&=&sdiam_k(P_{m_1}\Box P_{m_2}\Box \cdots \Box P_{m_{r-2}})+(k-2)(m_{r-1}-1)+(k-2)(m_{r}-1)\\
&\leq&\ldots\\
&\leq&sdiam_k(P_{m_1})+(k-2)\left(\sum_{i=2}^{r}m_i-r+1\right)\\
&=&(k-2)\left(\sum_{i=2}^{r}m_i-r+1\right)+m_1-1.
\end{eqnarray*}

$(2)$ From Theorem \ref{th2-2}, the result holds.\qed
\end{pf}

\subsection{$r$-dimensional torus}

An \emph{$r$-dimensional torus} is the Cartesian product of $r$
cycles $C_{m_1},C_{m_2},\cdots,C_{m_r}$ of size at least three. The
cycles $C_{m_i}$ are not necessary to have the same size.
\cite{Ku} showed that there are $r$ edge-disjoint spanning trees in
an $r$-dimensional torus. The network $C_{m_1}\circ C_{m_2}\circ
\cdots \circ C_{m_r}$ is investigated in \cite{Mao2}. Here, we
consider the networks constructed by $C_{m_1}\Box C_{m_2}\Box \cdots
\Box C_{m_r}$ and $C_{m_1}\circ C_{m_2}\circ \cdots \circ C_{m_r}$.
\begin{pro}\label{pro4-4}
Let $k,m_1,m_2,\cdots,m_r$ be the integers with $m_1\geq m_2 \geq
\cdots \geq m_r\geq 3$ and $3\leq k\leq \prod_{i=1}^rm_i$.

$(1)$ For network $C_{m_1}\Box C_{m_2}\Box \cdots \Box C_{m_r}$,
\begin{eqnarray*} \label{n-t=0}
\sum_{i=1}^{r}\left\lfloor\frac{(k-1)m_i}{k}\right\rfloor
&\leq&sdiam_k(C_{m_1}\Box C_{m_2}\Box \cdots \Box C_{m_r})\\
&\leq&\left \lfloor\frac{m_1(k-1)}{k}\right\rfloor
+(k-2)\sum_{i=2}^{n}\left \lfloor\frac{m_i(k-1)}{k}\right \rfloor,
\end{eqnarray*}
where $m_i$ is the order of $C_{m_i}$ and $1\leq i\leq n$.

$(2)$ For network $C_{m_1}\circ C_{m_2}\circ \cdots \circ C_{m_r}$,
$$
sdiam_k(C_{m_1}\circ C_{m_2}\circ \cdots \circ C_{m_r})\leq \left\{
\begin{array}{ll}
\left\lfloor\frac{(k-1)m_1}{k}\right\rfloor+k-2, &\mbox {\rm if} \ k\leq m_1;\\[0.2cm]
m_1+k-3, &\mbox {\rm if} \ m_1<k\leq \prod_{i=1}^rm_i.
\end{array}
\right.
$$
and
$$
sdiam_k(C_{m_1}\circ C_{m_2}\circ \cdots \circ C_{m_r})\geq \left\{
\begin{array}{ll}
\left\lfloor\frac{(k-1)m_1}{k}\right\rfloor, &\mbox {\rm if} \ 2\leq k\leq m_1;\\[0.2cm]
k-1, &\mbox {\rm if} \ m_1+1\leq k\leq \prod_{i=1}^rm_i.
\end{array}
\right.
$$
\end{pro}
\begin{pf}
$(1)$ From $(3)$ of Proposition \ref{pro4-1},
$sdiam_k(C_{m_i})=\left\lfloor\frac{(k-1)m_i}{k}\right\rfloor$ for
each $i \ (1\leq i\leq r)$. By Theorem \ref{th2-2}, we have
$$
sdiam_k(C_{m_1}\Box C_{m_2}\Box \cdots \Box C_{m_r})\geq
\sum_{i=1}^{r}sdiam_k(C_{m_i})=\sum_{i=1}^{r}\left\lfloor\frac{(k-1)m_i}{k}\right\rfloor.
$$
and
$$sdiam_k(C_{m_1}\Box C_{m_2}\Box \cdots \Box C_{m_r})
\leq
\left\lfloor\frac{(k-1)m_1}{k}\right\rfloor+(k-2)\sum_{i=2}^{r}\left\lfloor\frac{(k-1)m_i}{k}\right\rfloor.
$$

$(2)$ The result follows from Theorem \ref{th3-2}.\qed
\end{pf}

\subsection{$r$-dimensional generalized hypercube}

Let $K_m$ be a clique of $m$ vertices, $m\geq 2$. An
\emph{$r$-dimensional generalized hypercube} or \emph{Hamming graph}
\cite{DayA, Fragopoulou} is the product of $r$ cliques. We have the
following:
\begin{pro}\label{pro4-5}
Let $k,m_1,m_2,\cdots,m_r$ be the integers with $m_1\geq m_2 \geq
\cdots \geq m_r\geq k\geq 2$.

$(1)$ For network $K_{m_1}\Box K_{m_2}\Box \cdots \Box K_{m_r}$,
$$
r(k-1)\leq sdiam_k(K_{m_1}\Box K_{m_2}\Box \cdots \Box K_{m_r})\leq
(k-1)(kr-2r-k+3).
$$

$(2)$ For network $K_{m_1}\circ K_{m_2}\circ \cdots \circ K_{m_r}$,
$$
sdiam_k(K_{m_1}\circ K_{m_2}\circ \cdots \circ K_{m_r})=k-1.
$$
\end{pro}
\begin{pf}
$(1)$ From $(1)$ of Proposition \ref{pro4-1}, $sdiam_k(K_{m_i})=k-1$
for each $i \ (1\leq i\leq r)$. From Theorem \ref{th2-2}, we have
$$
sdiam_k(K_{m_1}\Box K_{m_2}\Box \cdots \Box K_{m_r})\geq
\sum_{i=1}^rsdiam_k(K_{m_i})=r(k-1)
$$
and
$$
sdiam_k(K_{m_1}\Box K_{m_2}\Box \cdots \Box K_{m_r})\leq
(k-2)(k-1)(r-1)+(k-1)=(k-1)(kr-2r-k+3).
$$

$(2)$ From the definition of lexicographical product, $K_{m_1}\circ
K_{m_2}\circ \cdots \circ K_{m_n}$ is a complete graph, and hence
$sdiam_k(K_{m_1}\circ K_{m_2}\circ \cdots \circ K_{m_n})=k-1$.\qed
\end{pf}

\subsection{$n$-dimensional hyper Petersen network}

An \emph{$n$-dimensional hyper Petersen network} $HP_n \ (n\geq 3)$
is defined as follows (see \cite{Das}).
\begin{itemize}
\item $HP_3$ is the Petersen graph (see Fig.4
$(a)$);
\begin{figure}[htbp]
  \begin{center}
    \includegraphics[scale=0.7]{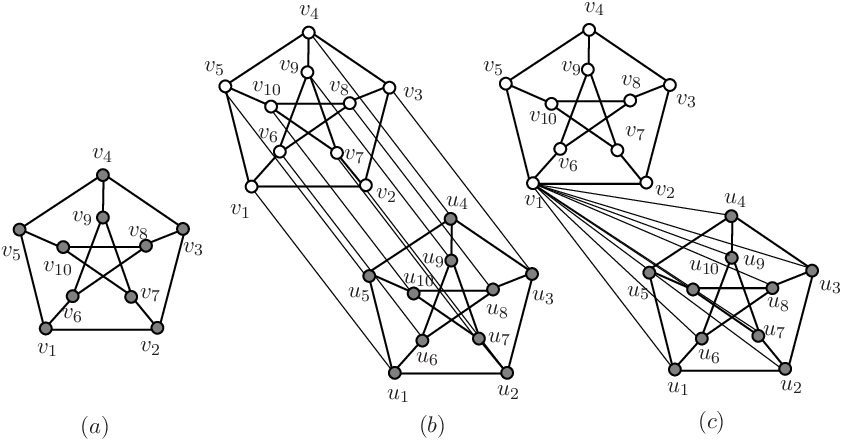}
    \caption{$(a)$ Petersen graph; $(b)$ The network $HP_4$; $(c)$ The
structure of $HL_4$.}
    \label{fig:logo}
  \end{center}
\end{figure}

\item $HP_n$ is the Cartesian product of the Petersen graph $PG$ and an
$(n-3)$-dimensional hypercube $Q_{n-3}$, that is, $HP_n=PG\Box
Q_{n-3}$, where $n\geq 4$.
\end{itemize}
The hyper Petersen network $HP_4$ are depicted in Fig.4 $(b)$.

The network $HL_n \ (n\geq 3)$ is defined as follows (see
\cite{Mao2}).
\begin{itemize}
\item $HL_3$ is the Petersen graph;

\item $HL_n$ is the lexicographic product of the Petersen graph $PG$ and an
$(n-3)$-dimensional hypercube $Q_{n-3}$, that is, $HP_n=PG\circ
Q_{n-3}$, where $n\geq 4$.
\end{itemize}
Note that $HL_4$ is a graph obtained from two copies of the Petersen
graph by add one edge between one vertex in a copy of the Petersen
graph and one vertex in another copy. See Figure 4 $(c)$ for an
example (We only show the edges $v_1u_i \ (1\leq i\leq 10)$).

Similarly to the proof of $(4)$ of Theorem \ref{th2-2}, we can get
the following observation.
\begin{obs}\label{obs4-1}
Let $G$ be a connected graph of order $n$. If $n-\kappa(G)+1\leq
k\leq n$, then $sdiam_k(G)=k-1$.
\end{obs}

\begin{pro}\label{pro4-6}
$(1)$ For network $HP_3$ and $HL_3$,
$$
sdiam_k(HP_3)=sdiam_k(HL_3)=\left\{
\begin{array}{ll}
k+1, &\mbox {\rm if} \ k=3,4;\\[0.1cm]
k, &\mbox {\rm if} \ k=5,6,7;\\[0.1cm]
k-1, &\mbox {\rm if} \ 8\leq k\leq 10.
\end{array}
\right.
$$

$(2)$ For network $HL_4$,
$$
sdiam_k(HL_4)=\left\{
\begin{array}{ll}
k, &\mbox {\rm if} \ 3\leq k\leq 7;\\[0.1cm]
k-1, &\mbox {\rm if} \ 8\leq k\leq 20.
\end{array}
\right.
$$

$(3)$ For network $HP_4$,
$$
\left\{
\begin{array}{ll}
sdiam_k(HP_4)=5, &\mbox {\rm if} \ k=3;\\[0.1cm]
k-1\leq sdiam_k(HP_4)\leq 9+\left\lfloor k/2\right\rfloor, &\mbox {\rm if} \ 4\leq k\leq 16;\\[0.1cm]
sdiam_k(HP_4)=k-1, &\mbox {\rm if} \ 17\leq k\leq 20.
\end{array}
\right.
$$
\end{pro}
\begin{pf}
$(1)$ Observe that $HL_3$ is just the Petersen graph. Set $G=HL_3$.
Choose $S=\{v_1,v_3,v_9\}$. One can see that any $S$-Steiner tree
must use at least $4$ edges of $G$, and hence $sdiam_3(G)\geq
d_G(S)\geq 4$. One can check that $d_G(S)\leq 4$ for any $S\subseteq
V(G)$ and $|S|=3$. Therefore, $sdiam_3(G)\leq 4$, and hence
$sdiam_3(G)=sdiam_3(HL_3)=4$. Since $HL_3=HP_3$, we have
$sdiam_3(HP_3)=sdiam_3(HL_3)=4$. Since $\kappa(G)=3$, it follows
from Observation \ref{obs4-1} that $sdiam_k(G)=k-1$ if $8\leq k\leq
10$. If $k=4$, then we choose $S=\{v_1,v_4,v_7,v_8\}$. One can see
that any $S$-Steiner tree must use at least $5$ edges of $G$, and
hence $sdiam_4(G)\geq d_G(S)\geq 5$. One can check that $d_G(S)\leq
5$ for any $S\subseteq V(G)$ and $|S|=3$. So, we have
$sdiam_4(G)=5$. Similarly, we can prove that $sdiam_k(G)=k$ if
$5\leq k\leq 7$.

$(2)$ For network $HL_4$, there are two copies of Petersen graphs,
say $HL_3$ and $HL'_3$. Set $G=HL_4$, $V(HL_3)=\{v_i\,|\,1\leq i\leq
10\}$ and $V(HL'_3)=\{u_i\,|\,1\leq i\leq 10\}$. Choose
$S=\{v_1,v_2,v_9\}$. One can see that any $S$-Steiner tree must use
at least $3$ edges of $G$, and hence $sdiam_3(G)\geq d_G(S)\geq 3$.
It suffices to show that $d_G(S)\leq 3$ for any $S\subseteq V(G)$
and $|S|=3$. Suppose $S\subseteq V(HL_3)$ or $S\subseteq V(HL'_3)$.
Without loss of generality, let $S=\{v_1,v_2,v_3\}\subseteq
V(HL_3)$. If $d_{HL_3}(S)=4$, then the tree induced by the edges in
$\{u_1v_1,u_1v_2,u_1v_3\}$ is an $S$-Steiner tree, and hence
$d_G(S)\leq 3$. Otherwise, $d_G(S)\leq d_{HL_3}(S)\leq 3$, as
desired. Suppose $|S\cap V(HL_3)|=2$ or $|S\cap V(HL'_3)|=2$.
Without loss of generality, let $|S\cap V(HL_3)|=2$ and
$S=\{v_1,v_2,u_1\}$. Then the tree induced by the edges in
$\{u_1v_1,u_1v_2\}$ is an $S$-Steiner tree, and hence $d_G(S)\leq
2$, as desired. So $sdiam_3(HL_4)=3$. Since $\kappa(G)=13$, it
follows from Observation \ref{obs4-1} that $sdiam_k(G)=k-1$ if
$8\leq k\leq 20$. One can also prove that $sdiam_k(G)=k$ if $3\leq
k\leq 7$.

$(3)$ For network $HP_4$, there are two copies of Petersen graphs,
say $HP_3$ and $HP'_3$. Set $G=HP_4$, $V(HP_3)=\{v_i\,|\,1\leq i\leq
10\}$ and $V(HP'_3)=\{u_i\,|\,1\leq i\leq 10\}$. Choose
$S=\{u_1,u_3,v_{10}\}$. One can see that any $S$-Steiner tree must
use at least $5$ edges of $G$, and hence $sdiam_3(G)\geq d_G(S)\geq
5$. One can check that $d_G(S)\leq 5$ for any $S\subseteq V(G)$ and
$|S|=3$. Then $sdiam_3(HP_4)\leq 5$, and hence $sdiam_3(HP_4)=5$.
Since $\kappa(G)=4$, it follows from Observation \ref{obs4-1} that
$sdiam_k(G)=k-1$ if $17\leq k\leq 20$. For $4\leq k\leq 16$, we have
$sdiam_k(HP_4)\geq k-1$, and for any $S\subseteq V(G)$ with $|S|=k$,
we let $S\cap V(HP_3)=S_1$ and $S\cap V(HP'_3)=S_2$. Without loss of
generality, let $|S_1|\geq \lceil\frac{k}{2}\rceil$. Let $S_2=S\cap
V(HP'_3)=\{u_1,u_2,\ldots,u_x\}$, where $x\leq \left\lfloor
k/2\right\rfloor$. Since $HP_3$ is connected, it follows that it
contains a spanning tree $T$ of size $9$. Then the tree induced by
the edges in $E(T)\cup \{u_iv_i\,|\,1\leq i\leq x\}$ is an
$S$-Steiner tree in $G$, and hence $d_G(S)\leq x+9\leq \left\lfloor
k/2\right\rfloor+9$. \qed
\end{pf}\\

\acknowledgements \label{sec:ack}
The authors are very grateful to
the referees' comments and suggestions, which helped to improve the
presentation of the paper.

\nocite{*}
\bibliographystyle{abbrvnat}

\begin{thebibliography}{49}
\providecommand{\natexlab}[1]{#1}
\providecommand{\url}[1]{\texttt{#1}}
\expandafter\ifx\csname urlstyle\endcsname\relax
  \providecommand{\doi}[1]{doi: #1}\else
  \providecommand{\doi}{doi: \begingroup \urlstyle{rm}\Url}\fi

\bibitem[Ali(2013)]{Ali}
P.~Ali.
\newblock The steiner diameter of a graph with prescribed girth.
\newblock \emph{Discrete Math.}, 313(12):\penalty0 1322--1326, 2013.

\bibitem[Ali et~al.(2012)Ali, Dankelmann, and Mukwembi]{AliDM}
P.~Ali, P.~Dankelmann, and S.~Mukwembi.
\newblock Upper bounds on the steiner diameter of a graph.
\newblock \emph{Discrete Appl. Math.}, 160(12):\penalty0 1845--1850, 2012.

\bibitem[Anand et~al.(2012)Anand, Changat, Klavar, and Peterin]{AnandCKP}
B.~Anand, M.~Changat, S.~Klavar, and I.~Peterin.
\newblock Convex sets in lexicographic product of graphs.
\newblock \emph{Graphs Combin.}, 28:\penalty0 77--84, 2012.

\bibitem[Bao et~al.(1998)Bao, Igarashi, and R.\"{O}hring]{Bao}
F.~Bao, Y.~Igarashi, and S.~R.\"{O}hring.
\newblock Reliable broadcasting in product network.
\newblock \emph{Discrete Appl. Math.}, 83:\penalty0 3--20, 1998.

\bibitem[Bondy and Murty(2008)]{Bondy}
J.~Bondy and U.~Murty.
\newblock \emph{Graph Theory}.
\newblock Springer, 2008.

\bibitem[Buckley and Harary(1990)]{Harary}
F.~Buckley and F.~Harary.
\newblock \emph{Distance in Graphs}.
\newblock Addision-Wesley, 1990.

\bibitem[C\'{a}ceresa et~al.(2008)C\'{a}ceresa, M\'{a}rquezb, and
  Puertasa]{Caceresa}
J.~C\'{a}ceresa, A.~M\'{a}rquezb, and M.~Puertasa.
\newblock Steiner distance and convexity in graphs.
\newblock \emph{European J. Combin.}, 29:\penalty0 726--736, 2008.

\bibitem[Calkin and Wilf(1998)]{Calkin}
N.~Calkin and H.~Wilf.
\newblock The number of independent sets in a grid graph.
\newblock \emph{SIAM J. Discrete Math.}, 11(1):\penalty0 54--60, 1998.

\bibitem[Chartrand et~al.(1989)Chartrand, Oellermann, Tian, and Zou]{Chartrand}
G.~Chartrand, O.~Oellermann, S.~Tian, and H.~Zou.
\newblock Steiner distance in graphs.
\newblock \emph{'{C}asopis pro p\v{e}stov\'{a}n\'{i} matematiky}, 114:\penalty0
  399--410, 1989.

\bibitem[Chartrand et~al.(2010)Chartrand, Okamoto, and Zhang]{ChartrandOZ}
G.~Chartrand, F.~Okamoto, and P.~Zhang.
\newblock Rainbow trees in graphs and generalized connectivity.
\newblock \emph{Networks}, 55:\penalty0 360--367, 2010.

\bibitem[Chung(1987)]{Chung}
F.~Chung.
\newblock Diameter of graphs: Old problems and new results.
\newblock In \emph{18th Southeastern Conf. on Combinatorics, Graph Theory and
  Computing}, 1987.

\bibitem[Dankelmann and Entringer(2000)]{DankelmannE}
P.~Dankelmann and R.~Entringer.
\newblock Average distance, minimum degree, and spanning trees.
\newblock \emph{J. Graph Theory}, 33:\penalty0 1--13, 2000.

\bibitem[Dankelmann et~al.(1999)Dankelmann, Swart, and
  Oellermann]{DankelmannSO2}
P.~Dankelmann, H.~Swart, and O.~Oellermann.
\newblock Bounds on the steiner diameter of a graph.
\newblock In \emph{Combinatorics and Graph Theory and Algorithms I and II},
  pages 269--279, Kalamazoo, MI, 1999. New Issues Press.

\bibitem[Das et~al.(1995)Das, \"{O}hring, and Banerjee]{Das}
S.~Das, S.~\"{O}hring, and A.~Banerjee.
\newblock into hyper petersen network: Yet another hypercube-like
  interconnection topology.
\newblock \emph{VLSI Design}, 2(4):\penalty0 335--351, 1995.

\bibitem[D'Atri and Moscarini(1988)]{Moscarini}
A.~D'Atri and M.~Moscarini.
\newblock Distance-hereditary graphs and, steiner trees and connected
  domination.
\newblock \emph{SIAM J. Comput.}, 17(3):\penalty0 521--538, 1988.

\bibitem[Day et~al.(1994)Day, Oellermann, and Swart]{DayOS}
D.~Day, O.~Oellermann, and H.~Swart.
\newblock Steiner distance-hereditary graphs.
\newblock \emph{SIAM J. Discrete Math.}, 7(3):\penalty0 437--442, 1994.

\bibitem[Day and Al-Ayyoub(1997)]{DayA}
K.~Day and A.-E. Al-Ayyoub.
\newblock The cross product of interconnection networks.
\newblock \emph{IEEE Trans. Parall. Distr. Sys.}, 8(2):\penalty0 109--118,
  1997.

\bibitem[Du et~al.(1993)Du, Lyuu, and Hsu]{Du}
D.~Du, Y.~Lyuu, and D.~Hsu.
\newblock Line digraph iteration and connectivity analysis of de {Bruijn} and
  {Kautz} graphs.
\newblock \emph{IEEE Trans. Comput.}, 42:\penalty0 612--616, 1993.

\bibitem[Fragopoulou et~al.(1996)Fragopoulou, Akl, and Meijer]{Fragopoulou}
P.~Fragopoulou, S.~Akl, and H.~Meijer.
\newblock Optimal communication primitives on the generalized hypercube
  network.
\newblock \emph{IEEE Trans. Parall. Distr. Comput.}, 32(2):\penalty0 173--187,
  1996.

\bibitem[Furtula et~al.(2016)Furtula, Gutman, and Katani\'{c}]{FurtulaGK}
B.~Furtula, I.~Gutman, and V.~Katani\'{c}.
\newblock Three-center harary index and its applications.
\newblock \emph{Iranian J. Math. Chem.}, 7(1):\penalty0 61--68, 2016.

\bibitem[Garey and Johnson(1979)]{GareyJ}
M.~Garey and D.~Johnson.
\newblock \emph{Computers and Intractibility: A Guide to the Theory of
  NP-Completeness}.
\newblock Freeman \& Company New York, 1979.

\bibitem[Goddard and Oellermann(2011)]{Goddard}
W.~Goddard and O.~Oellermann.
\newblock Distance in graphs.
\newblock In M.~Dehmer, editor, \emph{Structural Analysis of Complex Networks},
  pages 49--72. Birkh\"{a}user, Dordrecht, 2011.

\bibitem[Goddard and Oellrmann(1994)]{GoddardOS}
W.~Goddard and O.~Oellrmann.
\newblock Steiner distance stable graphs.
\newblock \emph{Discrete Math.}, 132:\penalty0 65--73, 1994.

\bibitem[Gologranc(2018)]{Gologranc}
T.~Gologranc.
\newblock Steiner convex sets and cartesian product.
\newblock \emph{Bull. Malays. Math. Sci. Soc.}, 41(2):\penalty0 627--636, 2018.

\bibitem[Gutman(2016)]{GutmanSDD}
I.~Gutman.
\newblock On steiner degree distance of trees.
\newblock \emph{Appl. Math. Comput.}, 283:\penalty0 163--167, 2016.

\bibitem[Gutman et~al.(2015)Gutman, Furtula, and Li]{GFL}
I.~Gutman, B.~Furtula, and X.~Li.
\newblock Multicenter wiener indices and their applications.
\newblock \emph{J. Serb. Chem. Soc.}, 80:\penalty0 1009--1017, 2015.

\bibitem[Hakimi(1971)]{Hakimi}
S.~Hakimi.
\newblock Steiner's problem in graph and its implications.
\newblock \emph{Networks}, 1:\penalty0 113--133, 1971.

\bibitem[Hammack et~al.(2011)Hammack, Imrich, and Klav\u{z}ar]{Hammack}
R.~Hammack, W.~Imrich, and S.~Klav\u{z}ar.
\newblock \emph{Handbook of product graphs}.
\newblock CRC Press, 2011.

\bibitem[Hsu(1994)]{Hsu}
D.~Hsu.
\newblock On container width and length in graphs, groups, and networks.
\newblock \emph{IEICE Transaction on Fundamentals of Electronics,
  Communications and Computer Science}, E77-A:\penalty0 668--680, 1994.

\bibitem[Hsu and {\L}uczak(1994)]{Hsu2}
D.~Hsu and T.~{\L}uczak.
\newblock Note on the $k$-diameter of $k$-regular $k$-connected graphs.
\newblock \emph{Discrete Math.}, 133:\penalty0 291--296, 1994.

\bibitem[Hwang et~al.(1992)Hwang, Richards, and Winter]{HwangRW}
F.~Hwang, D.~Richards, and P.~Winter.
\newblock \emph{The Steiner Tree Problem}.
\newblock North-Holland, Amsterdam, 1992.

\bibitem[Itai and Rodeh(1988)]{Itai}
A.~Itai and M.~Rodeh.
\newblock The multi-tree approach to reliability in distributed networks.
\newblock \emph{Inform. Comput.}, 79:\penalty0 43--59, 1988.

\bibitem[Johnsson and Ho(1989)]{Johnsson}
S.~Johnsson and C.~Ho.
\newblock Optimum broadcasting and personaized communication in hypercubes.
\newblock \emph{IEEE Trans. Comput.}, 38(9):\penalty0 1249--1268, 1989.

\bibitem[Ku et~al.(2003)Ku, Wang, and Hung]{Ku}
S.~Ku, B.~Wang, and T.~Hung.
\newblock Constructing edge-disjoint spanning trees in product networks.
\newblock \emph{IEEE Trans. Parall. Distr. Sys.}, 14(3):\penalty0 213--221,
  2003.

\bibitem[Levi(1971)]{Levi}
A.~Levi.
\newblock Algorithm for shortest connection of a group of graph vertices.
\newblock \emph{Sov. Math. Dokl.}, 12:\penalty0 1477--1481, 1971.

\bibitem[Li et~al.(2016)Li, Mao, and Gutman]{LMG}
X.~Li, Y.~Mao, and I.~Gutman.
\newblock The steiner wiener index of a graph.
\newblock \emph{Discuss. Math. Graph Theory}, 36:\penalty0 455--465, 2016.

\bibitem[Li et~al.(2017)Li, Mao, and Gutman]{LMG2}
X.~Li, Y.~Mao, and I.~Gutman.
\newblock Inverse problem on the steiner wiener index.
\newblock \emph{Discuss. Math. Graph Theory}, 38(1):\penalty0 83--95, 2017.

\bibitem[Mao()]{MaoSurvey}
Y.~Mao.
\newblock \emph{Steiner Distance in Graphs--A Survey}.
\newblock arXiv:1708.05779.

\bibitem[Mao(2016)]{Mao2}
Y.~Mao.
\newblock Path-connectivity of lexicographical product graphs.
\newblock \emph{Int. J. Comput. Math.}, 93(1):\penalty0 27--39, 2016.

\bibitem[Mao(2017)]{Mao}
Y.~Mao.
\newblock The steiner diameter of a graph.
\newblock \emph{Bull. Iran. Math. Soc.}, 43(2):\penalty0 439--454, 2017.

\bibitem[Mao and Das(2018)]{MaoDas}
Y.~Mao and K.~Das.
\newblock Steiner gutman index.
\newblock \emph{MATCH Commun. Math. Comput. Chem.}, 79(3):\penalty0 779--794,
  2018.

\bibitem[Mao et~al.(2016)Mao, Wang, and Gutman]{MWG}
Y.~Mao, Z.~Wang, and I.~Gutman.
\newblock Steiner wiener index of graph products.
\newblock \emph{Trans. Combin.}, 5(3):\penalty0 39--50, 2016.

\bibitem[Mao et~al.(2017{\natexlab{a}})Mao, Wang, Gutman, and
  Klobu\v{c}ar]{MWGK}
Y.~Mao, Z.~Wang, I.~Gutman, and A.~Klobu\v{c}ar.
\newblock Steiner degree distance.
\newblock \emph{MATCH Commun. Math. Comput. Chem.}, 78(1):\penalty0 221--230,
  2017{\natexlab{a}}.

\bibitem[Mao et~al.(2017{\natexlab{b}})Mao, Wang, Gutman, and Li]{MWGL}
Y.~Mao, Z.~Wang, I.~Gutman, and H.~Li.
\newblock Nordhaus-gaddum-type results for the steiner wiener index of graphs.
\newblock \emph{Discrete Appl. Math.}, 219:\penalty0 167--175,
  2017{\natexlab{b}}.

\bibitem[Mao et~al.(2018)Mao, Melekian, and Cheng]{MaoMC}
Y.~Mao, C.~Melekian, and E.~Cheng.
\newblock A note on the steiner $(n-k)$-diameter of a graph.
\newblock \emph{Int. J. Comput. Math.:CST}, 3(1):\penalty0 41--46, 2018.

\bibitem[Meyer and Pradhan(1987)]{Meyer}
F.~Meyer and D.~Pradhan.
\newblock Flip trees.
\newblock \emph{IEEE Trans. Comput.}, 37(3):\penalty0 472--478, 1987.

\bibitem[Oellermann(1995)]{Oellermann}
O.~Oellermann.
\newblock From steiner centers to steiner medians.
\newblock \emph{J. Graph Theory}, 20(2):\penalty0 113--122, 1995.

\bibitem[Oellermann(1999)]{Oellermann2}
O.~Oellermann.
\newblock On steiner centers and steiner medians of graphs.
\newblock \emph{Networks}, 34:\penalty0 258--263, 1999.

\bibitem[Oellermann and Tian(1990)]{OellermannT}
O.~Oellermann and S.~Tian.
\newblock Steiner centers in graphs.
\newblock \emph{J. Graph Theory}, 14(5):\penalty0 585--597, 1990.

\end{thebibliography}

\label{sec:biblio}

\end{document}